\documentclass[a4paper,11pt]{article}%,dvipdfm

\usepackage{mathrsfs} 

\usepackage{amsthm}
\usepackage{hyperref} 
\usepackage{amssymb}
\usepackage{latexsym}
\usepackage{amsmath}
\usepackage{amsfonts}
\usepackage{mathabx}
\usepackage[dvips]{graphicx}
\usepackage[utf8]{inputenc}

\usepackage[LGR,T1]{fontenc}%
\usepackage[francais,english]{babel}
\usepackage{url}
\usepackage{enumerate}
\usepackage{hyperref}
\usepackage{xcolor}
\usepackage{subfig}

\usepackage{fancyhdr}
\newtheorem{Theorem}{Theorem}[part]
\newtheorem{Definition}{Definition}[part]
\newtheorem{Proposition}{Proposition}[part]

\newtheorem{Assumption}{Assumption}[part]
\newtheorem{Lemma}{Lemma}[part]
\newtheorem{Corollary}{Corollary}[part]
\newtheorem{Remark}{Remark}[part]
\newtheorem{Example}{Example}[part]

\makeatletter \@addtoreset{equation}{section}

\@addtoreset{Definition}{section}

\@addtoreset{Theorem}{section}

\@addtoreset{Proposition}{section}

\@addtoreset{Property}{section}

\@addtoreset{Assumption}{section}

\@addtoreset{Corollary}{section}

\@addtoreset{Lemma}{section}

\@addtoreset{Remark}{section}

\@addtoreset{Example}{section}

\newcommand{\cA}{\mathcal{A}}

\newcommand{\cD}{\mathcal{D}}
\newcommand{\cE}{\mathcal{E}}
\newcommand{\cF}{\mathcal{F}}
\newcommand{\cG}{\mathcal{G}}
\newcommand{\cH}{\mathcal{H}}
\newcommand{\cI}{\mathcal{I}}

\newcommand{\cK}{\mathcal{K}}
\newcommand{\cL}{\mathcal{L}}
\newcommand{\cM}{\mathcal{M}}
\newcommand{\cN}{\mathcal{N}}

\newcommand{\cP}{\mathcal{P}}

\newcommand{\cS}{\mathcal{S}}
\newcommand{\cT}{\mathcal{T}}

\newcommand{\cV}{\mathcal{V}}
\newcommand{\cW}{\mathcal{W}}

\newcommand{\F}{\mathbb{F}}

\renewcommand{\P}{\mathbb{P}}

\newcommand{\R}{\mathbb{R}}

\def \proof{{\noindent \bf Proof. }}
\def \eproof{\hbox{ }\hfill$\Box$}

\newcommand{\ud}{\mathrm{d}}

\newcommand{\1}{{\bf 1}}
    
\newcommand{\set}[1]
    {\ensuremath{\{ #1 \}}}
    
\newcommand{\HP}[1] %L2DPDT sur 0,T
    {\ensuremath{\mathscr{H}^{#1}}}

\newcommand{\esp}[1]{\ensuremath{\mathbb{E} \!\! \left[#1\right] }}

\renewcommand{\Xi}[1]{X_{i #1}}

\title{Convergence of particles and tree based scheme for singular FBSDEs\\
}

\author{
  Jean-Fran{\c{c}}ois CHASSAGNEUX\thanks{UFR de Math{\'e}matiques \& LPSM, Universit{\'e} Paris Cit\'e, B{\^a}timent Sophie Germain, 8 place Aur{\'e}lie Nemours, 75013 Paris, France ({\tt chassagneux@lpsm.paris})} \footnotemark[3], Mohan YANG\thanks{UFR de Math{\'e}matiques \& LPSM, Universit{\'e} Paris Cit\'e, B{\^a}timent Sophie Germain, 8 place Aur{\'e}lie Nemours, 75013 Paris, France ({\tt myang@lpsm.paris})}  \footnote{The authors would like to thank Benjamin Jourdain for fruitful discussions.}}

\definecolor{r_one}{rgb}{0,0,156}
\definecolor{r_two}{HTML}{000000}

\date{}

\begin{document}
\maketitle

\begin{abstract}
We study an implementation of the theoretical splitting scheme introduced in \cite{chassagneux2022numerical} for singular FBSDEs \cite{carmona2013singular} and their associated quasi-linear degenerate PDEs. The fully implementable algorithm is based on particles approximation of the transport operator and tree like approximation of the diffusion operator appearing in the theoretical splitting. We prove the convergence with a rate of our numerical method under some reasonable conditions on the coefficients functions. This validates \emph{a posteriori} some numerical results obtained in \cite{chassagneux2022numerical}.
We conclude the paper with a numerical section presenting various implementations of the algorithm and discussing their efficiency in practice.
\end{abstract}
%\tableofcontents

\section{Introduction}
In this work, we study the theoretical and numerical convergence of a particles and tree based scheme for singular FBSDEs.
The class of singular Forward Backward Stochastic Differential Equations (FBSDE) is motivated by application in the modeling of carbon markets \cite{howison2012risk,carmona2013singularfirst,chassagneux2022modelling}. \textcolor{black}{In these models, the FBSDE solution represents the equilibrium price of an allowance (to emit one ton of $CO_2$) in a perfect market.} Let $(\Omega,\cF,\P)$ be a stochastic basis supporting a $d$-dimensional Brownian motion $W$ and $T>0$ a terminal time. We denote by $\mathbb{F}:=(\cF_t)_{t \ge 0}$ the filtration generated by the Brownian motion (augmented). The singular FBSDE system, with solution $(P_t,E_t,Y_t,Z_t)_{0\le t\le T}$, has the following form:
\begin{align} \label{eq singular fbsde}
	\left\{
	\begin{array}{rcl}
		dP_t & = & b(P_t)dt + \sigma(P_t) d W_{t}\\
		dE_t &=  &\mu(Y_t,P_{t})d t\\
		dY_t &=  & Z_t d W_t
	\end{array}
	\right. .
\end{align}
The terminal condition is typically given by 
\begin{align}\label{eq very special spec phi}
Y_T = \phi(E_T) \text{ where } \phi(e) = \1_{e \ge \Lambda}\,,
\end{align} 
for some $\Lambda > 0$.  Existence and uniqueness to the above FBSDE is not straightforward has it is fully coupled (the backward process $Y$ appears in the coefficient of the forward process $E$), it is degenerate in the forward direction $E$ and the terminal condition is discontinuous. However, assuming some structural conditions, see Theorem \ref{th existence uniqueness one-period}  below for a precise statement, Carmona and Delarue \cite{carmona2013singular} managed to prove the wellposedness of such singular FBSDEs. They also show that the following Markovian representation holds true for the $Y$-process, namely:
\begin{align}
Y_t := \cV(t,P_t,E_t) \;,\; 0 \le t < T,
\end{align}
where the measurable function
$\cV: [0,T]\times\R^d\times\R \rightarrow \R$
is classically named the \emph{decoupling field}.

%\vspace{2mm}
%{\color{magenta}
%Motivated by applications in the modeling of carbon markets, we will also consider the case of multi-period FBSDEs (for a finite number of period $Q$) Given $(T_q)_{0 \le q \le Q}$ with $T_0 =0$ and $T_q < T_{q+1}$, we consider on $[T_q, T_{q+1})$ that $Y$ satisfies \eqref{eq singular fbsde} and the transition from one period to another is given (backwardly) by
%\begin{align}
%Y_{T_q - } = Y_{T_q +}\1_{E_{T_q} \le \Lambda} + \1_{E_{T_q} \ge \Lambda}.
%\end{align}
%
%\textcolor{red}{add more multi-period case, recall that above equalities have to be understood in the weak sense}
%Note that in practice one has to work with ``regularised version'' of the PDE (In the proof below I will do as if we were using a regularised version of the PDE but this will have to be written properly). It could also be interesting to look at terminal condition 
%$(p,e) \mapsto \phi(p,e)$ to work with multi-period model. \\
%}

\noindent The decoupling field $\cV$ associated to \eqref{eq singular fbsde} is a weak solution to
%In this case, we denote $u$ the value function of the entropy solution to
\begin{align}\label{eq PDE general}
	\partial_t u + \partial_e ( \mathfrak{M}(p,u) )+ \partial_p u b(p) +\frac{1}{2}\textrm{Tr}[\sigma(p)\sigma^\top\!(p) \partial_{pp}^{2}u]  =0 \; \text{ and }\; u(T,p,e) = \phi(e),
\end{align}
with
\begin{align}
\mathfrak{M}(p,y)& :=\int^y_0 \mu(p,\upsilon) \ud \upsilon\;,\;  \label{eq de op conservative}
0\le y \le 1,
%\\
%\text{ and } \mathcal{L}^p &:= \langle b(p),\partial_p \rangle + \frac{1}{2}\textrm{Tr}[\sigma(p)\sigma^\top(p) \partial_{pp}^{2}]
%\label{eq de dynkin op}
\end{align}
 where \textcolor{black}{$\partial_p$ denotes the Jacobian with respect to $p$}, $\partial_t$ the time derivative, $\partial_e$ the derivative with respect to the $e$ variable,  $\top$ is the transpose and $\partial^{2}_{p p}$ is the matrix of second derivatives  with respect to the $p$ variable. 

The algorithm we study is based on a splitting scheme designed for this kind of equation and introduced in \cite{chassagneux2022numerical}. Though it could be defined at a general level for the PDE, we should note that the convergence proof relies on the FBSDE setting. Indeed, under appropriate conditions, the FBSDE appears as well-posed random characteristics for the PDE \eqref{eq PDE general}. The splitting scheme involves the iteration on a discrete time grid of the composition of two operators : a transport operator $\cT$ (in this case the $P$ variable is fixed) and a diffusion operator $\cD$ (in this case the $E$ variable is fixed). Precise definitions are given in the next section. From \cite{chassagneux2022numerical}, we already know that the theoretical splitting scheme is convergent with a rate one half. This result is obtained under the same structural conditions which guarantee the well-posedness of the FBSDE. In \cite{chassagneux2022numerical}, the splitting scheme is then implemented using various finite difference approximations of the transport operator and a non-linear regression using Deep Neural Networks for the diffusion operator.  An alternative scheme \cite[Section 3.3.1]{chassagneux2022numerical} is also suggested to verify the convergence result of the main numerical procedure. The goal of this paper is to analyse precisely this alternative scheme by proving its convergence with a rate under some regularity conditions.

\vspace{2mm}
The rest of the paper is organized as follows. In the next section, we recall the theoretical setting for the study of singular FBSDEs and recall the theoretical splitting scheme. \textcolor{black}{We also introduce the framework used to study the convergence of the numerical methods and give some properties of the solution $\cV$ in this framework.} In Section \ref{se algo}, we present the numerical algorithm and some of its variants that will be used for numerical tests. We also state our main convergence result in Theorem \ref{th main conv result}. Section \ref{se errors} is dedicated to the proof of the convergence by studying precisely all sources of errors.
The last section presents numerical results for the various schemes introduced, illustrating the theoretical convergence.

\paragraph{Notations}
In the  following we will  use the following spaces: 

\noindent Let $\mathbb{H}=(\cH_t)_{t \ge 0}$ denotes a generic right-continuous filtration,  
\begin{itemize}
\item 
%\begin{bluetext}
For  fixed $0\le a<b<+\infty$ and $I=[a,b]$ or $I=[a,b)$, { $\mathcal{S}^{2,k}(I,\mathbb{H})$ is the
set of  $\mathop{}\!\mathbb{R}^k$-valued c\`adl\`ag\footnote{French acronym for right continuous with left limits.} $\mathcal{H}_t$-adapted processes $Y$, s.t.\\[-4mm]  
\begin{align*}
\|Y\|_{\mathcal{S}^2}^2:= \mathbb{E}\left[{\sup_{t \in I} |Y_t|^2 } %
\right] <\infty .
\end{align*}
Note that we may omit the dimension and the terminal date in the norm notation as this will be clear from the context. $\mathcal{S}^{2,k}_\mathrm{c}(I,\mathbb{H})$ is the {subspace of processes} with continuous sample paths.} 
%\\
%We also consider ${\mathcal{S}}^{2,k}( [0,\infty))$ the vector space of c\`adl\`ag adapted processes $Y$, with values in $\R^k$, and such that $\mathbb{E}\left[\sup_{0 \leq t  \leq b} |Y_t|^2  \right] <\infty $ for every $b > 0$.  ${\mathcal{S}}^{2,k}_{\mathrm{c}}( [0,\infty))$ denotes the subspace of such processes having continuous paths.

\item For fixed $0\le a<b<+\infty$, and $I=[a,b]$, we denote by { $\mathcal{H}^{2,k}(I,\mathbb{H})$ the set
of $\mathop{}\!\mathbb{R}^k$-valued progressively measurable processes $Z$, such that \\[-4mm] 
\begin{equation*}
\|Z\|_{\mathcal{H}^2}^2 := \mathbb{E} \left[{\int_I |Z_t|^2 dt} \right] 
<\infty.  
\end{equation*}%
%$\mathcal{H}^{2,k}([0,\infty))$ is the set
%of $\mathop{}\!\mathbb{R}^k$-valued progressively measurable processes $Z$, such that 
%$ \mathbb{E} \left[{\int_0^b |Z_t|^2 dt} \right] 
%<\infty$, for all $b>0$.
}
%\end{bluetext}
\end{itemize}
%\noindent Note that we will omit to specify $\mathbb{H}$ when 

\noindent We denote by $\mathscr{P}(\R)$ the space of probability measure on $\R$ and $\mathscr{P}_q(\R)$, $q\ge 1$ the subset of probability measure which have a $q$-moment finite. 
%These spaces are endowed with the corresponding Wasserstein distance denoted $\cW_q$.
\\
We denote by $\mathscr{I}$ the space of cumulative distribution function (CDF) namely functions $\theta$ given by
$$\R \ni x \mapsto \theta(x) = \mu((-\infty,x]) \in [0,1]$$
for some $\mu \in \mathscr{P}(\R)$.
\\
{\color{black} The set $\mathscr{C}^{1,2,1}(I\times\R^d \times \R)$ (resp. $\mathscr{C}^{1,2,2}(I\times\R^d \times \R)$) with $I=[0,T)$ or $I = [0,T]$, is the set of functions $f:I\times\R^d\times \R \rightarrow \R$ continuously differentiable: once in their first variable, twice in their second variable and once (resp. twice) in their third variable. 
}

\section{Review of theoretical result for singular FBSDEs}

In this section, we first recall  the main properties of solution to \eqref{eq singular fbsde}. We then define the theoretical splitting scheme studied in \cite{chassagneux2022numerical} and used to design our implemented algorithms introduced in the next section. Finally, we consider a smooth framework for the solution  of \eqref{eq PDE general} and obtain some properties useful to prove the convergence of the numerical schemes. 

\subsection{Well-posedness of Singular FBSDEs}
We first define a class of function to which the solution belongs.
\begin{Definition}\label{de cK}
Let $\cK$ be the class of functions $\Phi:\R^d\times\R \rightarrow [0,1]$ such that $\Phi$ is $L_\phi$-Lipschitz  in the first variable for some $L_\Phi>0$, namely
\begin{align}
|\Phi(p,e)-\Phi(p',e)| &\le L_\Phi|p-p'| \quad\text{ for all }\quad (p,p',e) \in \R^d\times\R^d\times\R\;,
%\\
%\phi(p,e') &\ge \phi(p,e) \quad\text{ if }\quad e' \ge e\;,
\end{align} 
and moreover  satisfying, for all $p \in \R^d$,
\begin{align}
\Phi(p,e) = \nu(p,(-\infty,e]) \text{ where } \nu(p,\cdot) \in \cP_0(\R)\;.
\end{align}
%
%\textcolor{black}{
%\begin{align}\label{eq condition on phi}
%\sup_e \phi(p,e) = 1  \;\text{ and }\; \inf_e \phi(p,e) = 0 \quad\text{ for all }\quad p \in \R^d\;.
%\end{align}
%Do we need that ?=> $\phi$ bounded should be OK? no need to specify these bounds that exists anyway as $\phi$ is increasing...
%}
%\textcolor{red}{add right continuous or define directly as a cdf w.r.t $e$?}
\end{Definition}
\noindent
From the above definition, we have that $\Phi(p,\cdot) \in \mathscr{I}$ and is thus a non-decreasing, right continuous function and satisfies, for all $p \in \R^d$,
\begin{align} \label{eq some properties of Phi}
\lim_{e \rightarrow -\infty} \Phi(p,e) = 0 \text{ and } \lim_{e \rightarrow -\infty} \Phi(p,e) = 1\,.
\end{align}

\noindent We now recall the existence and uniqueness result for the singular FBSDE. It is also key to define the theoretical splitting scheme. This result is obtained for the following class of admissible coefficient functions.

\begin{Definition} \label{de class coef} Let $\cA$ be the class of functions 
$B:\R^d \rightarrow \R^d$, $\Sigma:\R^d \rightarrow \cM_d$, $F:\R\times \R^d\rightarrow \R$ which are  $L$-Lipschitz continuous functions. Moreover, \textcolor{black}{$F$ is strictly decreasing in $y$} and satisfies, for all $p \in \R^d$,
\begin{align}\label{eq conservation law}
\ell_1 |y-y'|^2 \le (y-y')(F(y',p)-F(y,p)) \le \ell_2 |y-y'|^2,
\end{align}
where  $L, \ell_1$ and $\ell_2$ are positive constants.
\end{Definition}

\noindent The well-posedness result in this setting reads as follows.
%\begin{magetext}
\begin{Theorem}[Proposition 2.10 in \cite{carmona2013singular}, Proposition 3.2 in \cite{chassagneux2022modelling}] \label{th existence uniqueness one-period} 
Let $\tau > 0$, $(B,\Sigma,F) \in \cA$ and $\Phi \in \cK$.
%:\R^d \times \R \rightarrow \R$ be a bounded measurable function, non-decreasing in its second variable, and  satisfying,
%\begin{align*}
%|\Phi(p,e)-\Phi(p',e)| \le L |p-p'|\;,\; \text{ for all } (p,p',e) \in \R^d\times\R^d \times \R\,.
%\end{align*} 

\noindent Given any initial condition $(t_0,p,e) \in [0,\tau) \times \mathbb{R}^{d} \times \mathbb{R}$, there exists a unique progressively measurable 4-tuple of processes $(P^{t_0,p,e}_t, E^{t_0,p,e}_t, Y^{t_0,p,e}_t, Z^{t_0,p,e}_t)_{t_0 \leq t \leq \tau} \in \mathcal{S}^{2,d}_{\mathrm{c}}([t_0,\tau],\F) \times \mathcal{S}^{2,1}_{\mathrm{c}}([t_0,\tau],\F) \times \mathcal{S}^{2,1}_{\mathrm{c}}([t_0,\tau),\F) \times \mathcal{H}^{2,d}([t_0,\tau],\F)$ satisfying the dynamics
\begin{align}  \label{general single period fbsde with p}
\begin{aligned}
\ud P^{t_0,p,e}_t &= B(P^{t_0,p,e}_t) \ud t + \Sigma(P^{t_0,p,e}_t) \ud W_t, & P^{t_0,p,e}_{t_0} &= p \in \mathbb{R}^{d}, \\
\ud E^{t_0,p,e}_t &= F(P^{t_0,p,e}_t, Y^{t_0,p,e}_t) \ud  t, & E^{t_0,p,e}_{t_0}
&= e \in \mathbb{R}, \\ 
\ud Y^{t_0,p,e}_t &=  Z^{t_0,p,e}_{t}   \ud {W_t}, & &
\end{aligned}
\end{align}
and such that
	\begin{align}
	\label{relaxed terminal condition with p}
	\mathbb{P} \left[ \Phi_{-}(P^{t_0,p,e}_\tau,E^{t_0,p,e}_\tau) \leq \lim_{t \uparrow \tau}Y^{t_0,p,e}_t \leq \Phi_{+}(P^{t_0,p,e}_\tau, E^{t_0,p,e}_\tau) \right] = 1.
	\end{align}
The unique \emph{decoupling field} defined by
\begin{align*}
[0,\tau) \times \mathbb{R}^{d} \times \mathbb{R} \ni (t_0,p,e) \rightarrow w(t_0,p,e) = Y^{t_0,p,e}_{t_0} \in \R
\end{align*}	
is continuous and satisfies
\begin{enumerate}
		\item{\label{v lipschitz constant e}For any $t \in [0,\tau)$, the function $w(t,\cdot, \cdot)$ is $1/(l_1(\tau-t))$-Lipschitz continuous with respect to $e$,}
		\item{\label{v lipschitz constant p}For any $t \in [0,\tau)$, the function $w(t,\cdot, \cdot)$ is $C$-Lipschitz continuous with respect to $p$, where $C$ is a constant depending on $L$, $\tau$ and $L_{\phi}$ only.}
		\item{\label{v condition terminal} Given $(p,e) \in \R^d \times \R$, for any family $(p_t, e_t)_{0 \leq t < \tau}$ converging to $(p,e)$ as $t \uparrow \tau$, we have
		\begin{align}\label{eq v condition terminal}
		\Phi_-(p,e) \le \liminf_{t \rightarrow \tau}w(t,p_t,e_t) \le \limsup_{t \rightarrow \tau}w(t,p_t,e_t) \le \Phi_+(p,e)\,.
		\end{align}
		\item{\label{v belongs to K}}
		\textcolor{black}{For any $t \in [0,\tau)$, the function $w(t,\cdot, \cdot) \in \cK$.}
}
\end{enumerate}
%\textcolor{red}{For all $t \in [0,T)$, $v(t,\cdot)\in \cK$.}
%\textcolor{magenta}{For later use, we define $\Theta_\tau (\Phi) = v(0,\cdot)$ with the parameters $B := b$, $\Sigma := \sigma$ and $F = -\mu$. This is the semi-group associated to the true solution.}
\end{Theorem}

%\end{magetext}

\noindent The previous result leads us to define a non-linear operator associated to \eqref{eq singular fbsde} under the following assumption, that will hold true for the rest of the paper:

\vspace{2mm}
%\noindent From now on, we make the following\\
\noindent \textbf{Standing assumptions:} The coefficients $(b,\sigma,\mu)$ of the singular FBSDE \eqref{eq singular fbsde} belong to $\cA$.

\begin{Definition}\label{de ope true solution}
We define the operator $\Theta$ by
\begin{align}
(0,{\infty})\times \cK \ni (h,\psi) \mapsto \Theta_h(\psi) = v(0,\cdot) \in \cK\,
\end{align}
where $\upsilon$ is the decoupling field given in Theorem \ref{th existence uniqueness one-period} with parameters $\tau = h$, $B=b$, $\Sigma = \sigma$, $F=\mu$ and $\Phi = \psi$.
\end{Definition}
%\noindent We also deduce from Theorem \ref{th existence uniqueness one-period}  that $(\Theta_t)_{0<t}$ is a semi-group of %non-linear operators. 

\noindent Let us now be given a discrete time grid of $[0,T]$:
\begin{align}
\pi := \set{t_0 := 0 \le \dots\le t_n\le \dots t_N := T}\,,
\end{align}
for $N\ge 1$. %We denote $|\pi| := \max_{n < N} (t_{n+1} - t_n)$.\\
For ease of presentation, we assume that the time grid $\pi$ is equidistant and thus, for $0\le n \le N-1$,
\begin{align}
t_{n+1} - t_n =\frac{T}N=:\mathfrak{h}\,.
\end{align}

\noindent Using Definition \ref{de ope true solution}, we observe that 
\begin{align}\label{eq V as iter nonlin op}
\cV(0,\cdot) := \Theta_T(\phi) = \prod_{0 \le n  < N}\Theta_{t_{n+1}-t_n}(\phi).
\end{align}
%, recall \eqref{eq de decoupling field}.

\vspace{2mm}
\noindent We recall now a result that arises from the proof of the previous Theorem \ref{th existence uniqueness one-period}, see \cite{chassagneux2022modelling}.

\begin{Corollary}[Approximation result]\label{co smooth approx of v}
Let $\tau > 0$, $(B,\Sigma,F) \in \cA$ and $\Phi \in \cK$. Let $(\Phi^k)_{k \ge 0}$ be a sequence of smooth functions belonging to $\cK$ and converging pointwise towards $\Phi$ as k goes to $+\infty$. For $\eta>0$, consider then $w^{\eta,k}$ the solution to:
\begin{align}
\label{value function pde approx}
\partial_t u +{F(u,p)} \partial_{e} u + \cL_p u + \frac12\eta^2 (\partial^2_{ee}u+ \Delta_{pp}u) = 0 \; \text{ and } \; u(\tau,\cdot) = \Phi^k
\end{align}
where $\Delta_{pp}$ is the Laplacian with respect to $p$, and $\cL_p$ is the operator
\begin{align}
\cL_p(\varphi)(t,p,e) = \partial_p\varphi (t,p,e)B(p) + \frac{1}{2} \mathrm{Tr} \left[ A(p) \partial^{2}_{pp} \right](\varphi)(t,p,e),
\end{align}
with $A = \Sigma \Sigma^{\top}$. For later use, we also define
$ 
\cL^\eta := \cL_p  + \frac12\eta^2 (\partial^2_{ee}+ \Delta_{pp})\;.
$
Then the functions $w^{\eta,k}$ are $\mathscr{C}^{1,2,2}([0,\tau]\times\R^d \times \R)$ and $\lim_{k \rightarrow \infty} \lim_{\eta \rightarrow 0} w^{\eta,k} = w$ where the convergence is locally uniform in $[0,\tau)\times\R^d \times \R$.
%\textcolor{red}{
Moreover, for all $k,\eta$, $w^{\eta,k}(t,\cdot) \in \cK$.
%}
\end{Corollary}

%{\color{magenta}
\noindent For later use, we introduce the system of FBSDEs associated to \eqref{value function pde approx}. We thus consider $(W,\tilde{B},\tilde{\tilde{B}})$ to be independent Brownian motions and denote $\mathbb{G}=(\cG_t)_{t\ge 0}$ the filtration they generate (augmented). For $t_0 \in [0,T)$ and $(p,e) \in \R^d\times \R$, let 
$$(P^\eta,E^{\eta,k}) \in \mathcal{S}^{2,d}_{\mathrm{c}}([t_0,T],\mathbb{G}) \times \mathcal{S}^{2,1}_{\mathrm{c}}([t_0,T],\mathbb{G}) $$ and 
$$(Y^{\eta,k},Z^{\eta,k},\tilde{Z}^{\eta,k},\tilde{\tilde{Z}}^{\eta,k}) \in \mathcal{S}^{2,d}_{\mathrm{c}}([t_0,T],\mathbb{G})  \times \mathcal{H}^{2,d}([t_0,T],\mathbb{G}) \times \mathcal{H}^{2,d}([t_0,T],\mathbb{G}) \times \mathcal{H}^{2,1}([t_0,T],\mathbb{G})$$ be the solution to
\begin{equation}\label{eq smoothed fbsde}
\left \{
\begin{split}
\ud P^\eta_t & = B(P^\eta_t) \ud t + \Sigma(P^\eta_t)\ud W_t + \eta \ud {\tilde{B}}_t,\;P^\eta_{t_0} = p \,,
\\
\ud E^{\eta,k}_t & = F(P^\eta_t,Y^{\eta,k}_t)\ud t + \eta \ud \tilde{\tilde{B}}_t,\; E^\eta_{t_0} = e \,,
\\
\ud Y^{\eta,k} &= Z^{\eta,k} \ud W_t + \tilde{Z}^{\eta,k}\ud {\tilde{B}}_t + \tilde{\tilde{Z}}^{\eta,k} \tilde{\tilde{B}}_t.
\end{split}
\right.
\end{equation}
Observe that $Y^{\eta,k}_t = w^{\eta,k}(t,P^\eta_t, E^{\eta,k}_t)$ for $t \in [t_0,T]$, in particular $Y^{\eta,k}_T = \Phi^k(P_T,E_t) $. Let us note that the above processes depend upon the initial condition $(t_0,p,e)$ but we  slightly abuse the notation in this regard by omitting to indicate the initial condition.
%}%ec

%\textcolor{red}{definition processus tangent lorsque les coeffs son assez réguliers ? si besoin plusieurs fois sinon dns preuve}

\subsection{A theoretical splitting scheme}
The numerical algorithm is based on a splitting method for the quasilinear PDE \eqref{eq PDE general}, which has been introduced in \cite{chassagneux2022numerical}. We now recall the splitting approach at a theoretical level.

%Using the previous result, we define the following operator associated to \eqref{eq singular fbsde}.
%\begin{Definition}\label{de ope true solution}
%We define the operator $\Theta$ by
%\begin{align}
%(0,{\infty})\times \cK \ni (h,\psi) \mapsto \Theta_h(\psi) = v(0,\cdot) \in \cK\,
%\end{align}
%where $\upsilon$ is the decoupling field given in Theorem \ref{th existence uniqueness one-period} with parameters $\tau = h$, $B=b$, $\Sigma = \sigma$, $F=\mu$ and $\Phi = \psi$.
%\end{Definition}
%We also deduce from Theorem \ref{th existence uniqueness one-period}  that $(\Theta_t)_{0<t}$ is a semi-group of non-linear operators. 
%In particular, we observe that $\cV(0,\cdot) := \Theta_T(\phi) = \prod_{0 \le n  < N}\Theta_{t_{n+1}-t_n}(\phi)$, recall \eqref{eq de decoupling field}.

\vspace{2mm}
\noindent We first introduce the transport step where the diffusion part is frozen.

\begin{Definition}[Transport step] \label{de transport step}
 %For $h \le T$ and $\psi \in \cK$, 
 We set
 $$(0,\infty) \times \cK \ni (h,\psi) \mapsto \cT_h(\psi) = \tilde{v}(0,\cdot) \in \cK\,,$$ 
 where $\tilde{v}$ is the decoupling field defined in Theorem \ref{th existence uniqueness one-period} with parameters
 $\tau = h$, $B =0$, $\Sigma = 0$, $F = \mu$ and terminal condition $\Phi = \psi$.
\end{Definition}

%{\color{blue}
\noindent In the definition above, $\tilde{v}(\cdot)$, defined on $[0,h]\times \R^d \times \R$, is the unique entropy solution to 
\begin{align} \label{eq de cons law associated with transport}
\partial_t w + \partial_e ( \mathfrak{M}(p,w) ) = 0\; \text{ and } \;\tilde{v}(h,\cdot)=\psi\,,
\end{align}
recall \eqref{eq PDE general}. 
\textcolor{black}{
\begin{Remark}\label{re abuse transport op} For later use in the convergence analysis, we introduce a small modification of the previous operator that acts on function in $\mathscr{I}$ instead of $\cK$. Namely, we set
 $$(0,\infty) \times \R^d \times \mathscr{I} \ni (h,p,\theta) \mapsto \tilde{\cT}_h(p,\theta) = \hat{v}_p(0,\cdot) \in \mathscr{I} $$ 
where $\hat{v}_p(0,\cdot)$ is 
the unique entropy solution to 
\begin{align} \label{eq de cons law associated with transport bis}
\partial_t w + \partial_e ( \mathfrak{M}(p,w) ) = 0\; \text{ and } \;\hat{v}_p(h,\cdot)=\theta.
\end{align}
One observes that $ \cT_h(\psi)(p,\cdot) = \tilde{\cT}_h(p,\psi(p,\cdot))$.
%the above is well defined from \eqref{eq de cons law associated with transport} where $p$ is indeed fixed.
\end{Remark}
}

\noindent Let us mention that, from e.g. \cite[Theorem 1.1]{jourdainreygner16}, we have the following stability result:
\begin{align}\label{eq L1 stab transport}
\int|\tilde{\cT}_h(p,\psi_1(p,\cdot))(e) - \tilde{\cT}_h(p,\psi_2(p,\cdot))(e)|\ud e
\le
\int |\psi_1(p,e) ) -\psi_2(p,e) | \ud e\,,
\end{align}
for $\psi_1,\psi_2 \in \cK$.

%We will use this fact in the numerical section.
%}

%\noindent Note, that, in particular, the following characteristic is well defined:
%\begin{align}
%\tilde{E}^{p,e}_s = e - \int_0^s f(\tilde{V}^{p,e}_s,P_s)\ud s\,,
%\end{align}
%where $\tilde{V}^{p,e}_s = \tilde{V}^{p,e}_s=\tilde{v}(s,p,\tilde{E}^{p,e}_s)$.

\vspace{2mm} 
\noindent We now introduce the diffusion step, where conversely, the $E$ - process is frozen to its initial value.
\begin{Definition}[Diffusion step]\label{de diffusion step}
 We set
 $$(0,\infty) \times \cK \ni (h,\psi) \mapsto \cD_h(\psi) = \bar{v}(0,\cdot) \in \cK$$ 
 where $\bar{v}(0,\cdot)$ is the decoupling field in Theorem \ref{th existence uniqueness one-period} with parameters
 $\tau = h$, $B =b$, $\Sigma = \sigma$, $F = 0$ and terminal condition $\Phi = \psi$.
\end{Definition}
\noindent Observe that, for $t \in [0,h)$,
\begin{align}\label{eq other de bar v}
\bar{v}(t,p,e) = \esp{\psi(P_h^{t,p},e)}\,,\; (p,e) \in \R^d \times \R \,\text{ and }  \bar{v}(t,\cdot) \in \cK\;.
\end{align}

\vspace{2mm}
\noindent 
%\textcolor{magenta}{
%Let us now be given a discrete time grid of $[0,T]$:
%\begin{align}
%\pi := \set{t_0 := 0 \le \dots\le t_n\le \dots t_N := T}\,,
%\end{align}
%for $N\ge 1$. %We denote $|\pi| := \max_{n < N} (t_{n+1} - t_n)$.\\
%For ease of presentation, we assume that the time grid $\pi$ is equidistant and thus, for $0\le n \le N-1$,
%\begin{align}
%t_{n+1} - t_n =\frac{T}N=:\mathfrak{h}\,.
%\end{align}
%}
\vspace{2mm}
\noindent We can now define the theoretical scheme on $\pi$ by a backward induction. Since
\begin{align}
\cV(0,\cdot) = \Theta_T(\phi) = \prod_{0 \le n  < N}\Theta_{t_{n+1}-t_n}(\phi)\,,
\end{align}
the main point of the scheme is to replace one step of $\Theta$ by one step of $\cT \circ \cD$. This leads to the following.

\begin{Definition}[Theoretical splitting scheme] \label{de th splitting theo}
We set 
\begin{align*}
(0,\infty)\times \cK \ni (h,\psi) \mapsto \cS_h(\psi) := \cT_h \circ \cD_h(\psi) \in \cK. 
\end{align*}

\noindent For $n \le N$, we denote by $u^\pi_n$ the solution of the following backward induction on $\pi$:
\\
- for $n=N$, set $u^\pi_N := \phi$,
\\
- for $n < N$, $u^\pi_n = \cS_{t_{n+1}-t_n}(u^\pi_{n+1})$.
\end{Definition}
\noindent The $(u^\pi_n)_{0 \le n \le N}$ stands for the approximation of the decoupling field $\cV(t,\cdot)$ for $t \in \pi$. 

%{\color{yellow}
%Moreover, we observe, from the property of $\cT$ and $\cD$, that 
%\begin{align}\label{eq u in K}
%u^\pi_n \in \cK, \text{ for all } 0 \le n \le N.
%\end{align}
%}

\noindent It is proven in \cite[Proposition 2.1]{chassagneux2022numerical} a key result concerning the scheme's truncation error.
\begin{Proposition}[truncation error]\label{pr truncation error}
Under our {standing assumptions on the coefficients} $(\mu,b,\sigma)$, the following holds, for  $\psi \in \cK$:
% Set $\tilde{v}(0,.) = \cS_h(\psi)$ and ${v}(0,.) = \Theta_h(\psi)$, then the following holds 
\begin{align}\label{eq pr truncation}
\int |\cS_h(\psi)(p,e)-\Theta_h(\psi)(p,e)| \ud e \le C_{L_\psi}(1+|p|^2)h^{\frac32}\;,
\end{align}
for $p \in \R^d$, $h>0$.
%\textcolor{magenta}{ $\tilde{C}$  depend on the Lipschitz constant of $\psi$ and $C$.}
\end{Proposition}
\noindent This allows in particular to obtain the convergence with a rate $\frac12$ of the theoretical splitting scheme to the decoupling field $\cV$, namely
\begin{align*}
\int_\R|\cV(0,p,e)-u^\pi_0(p,e)|\ud e \le C(1+|p|^2) \sqrt{\mathfrak{h}}\,,\quad \;p \in \R^d\,,
\end{align*}
for a positive constant $C$, see \cite[Theorem 2.2]{chassagneux2022numerical}.%, $p \in \R^D$.
%locally uniformly in $p$.
%\end{Theorem}

\subsection{Smooth setting for convergence results}

We now consider a more restricted framework.
In this setting, we give further properties of the value function $\cV$  that will be used to prove our convergence results.
Indeed, even if the scheme can be defined in a quite general setting, the convergence with a rate is obtained in a more regular framework.

\vspace{2mm}
\noindent Therefore, the main assumption we shall use is the following, see Remark \ref{re example for ass} below.

\begin{Assumption}\label{ass restrict regul} 
The \emph{decoupling field} $\cV$ is a $\mathscr{C}^{1,2,1}([0,T)\times\R^d \times \R)$ solution to
\begin{align}\label{eq smooth version}
\partial_t \cV +{\mu(\cV,p)} \partial_{e} \cV + b(p)\partial_p \cV+ \frac12\bar{\sigma}^2\Delta_{p}\cV = 0\;.
\end{align}
In particular, the function  $\sigma(\cdot)$ is constant $\sigma(\cdot) := \bar{\sigma}I_d$ with $\bar{\sigma}>0$ and $I_d$ the $d \times d$ identity matrix, recall \eqref{eq PDE general}.  
%namely 
%\begin{align}
%\cL^p u = b(p)\partial_p u + \frac12\bar{\sigma}^2 \Delta_{p}u \;.
%\end{align}
Moreover, the function $b$ is twice differentiable  with bounded and Lipschitz derivatives and the function $\mu$ is differentiable with Lipschitz and  bounded derivatives. 
\end{Assumption}

We first list some further properties of the function $\cV$ in the regular setting of the previous assumption.
\begin{Proposition}\label{pr property of V smooth}
Under Assumption \ref{ass restrict regul}, the function $\cV$ satisfies, for $(t,p,e) \in [0,T)\times\R^d\times\R$,
\begin{align}\label{eq control deriv}
0 \le \partial_e \cV(t,p,e) \le \frac{C}{T-t}\quad \text{ and } \quad
%\\
% &
 |\partial_p \cV(t,p,e)|  \le C\;.
\end{align}
Moreover, for $(t,p) \in [0,T)\times \R^d$,
\begin{align}\label{eq control deriv L1}
\int |\partial_p \cV(t,p,e)| \ud e \le C(T-t)\;,
\end{align}
and, thus, for $p' \in \R^d$,
\begin{align}\label{eq lipschitz L1}
\int |\cV(t,p,e)-\cV(t,p',e)| \ud e \le C(T-t)|p-p'| .
\end{align}
%\textcolor{red}{pour la seconde inégalité c'est multiplié par $T-t$ car $\partial_p\phi =0$ mais ça ne va pas vraiment servir...}
\end{Proposition}

%\begin{bluetext}
\proof The estimates \eqref{eq control deriv} are obtained directly from the property of $\cV$ given in Theorem \ref{th existence uniqueness one-period}. We now study the $L^1$-control.\\
Let $\cV^{\eta,k}$ be the smooth approximation given in Corollary \ref{co smooth approx of v} to $\cV$. We consider also the associated FBSDE $(P^\eta,E^{\eta,k},Y^{\eta,k})$ in \eqref{eq smoothed fbsde}. We adapt the computations in the proof of Lemma 3.8 in \cite{chassagneux2022modelling}, to obtain the equation corresponding to equation (3.35) there. In our simpler framework ($r=0$, $\phi$ does not depend on $p$), it reads
% in  \cite{jeanfrancois2020modelling} with $U_T = 0$ and $r=0$ and observing that $\cE = \partial_e E$, we obtain for any $0\le t_0\le T$,
\begin{align} \label{eq starting point}
\partial_p \cV^{\eta,k} (t_0,p,e) = \esp{\int_{t_0}^T\partial_p\mu(P^\eta_t,Y^{\eta,k}_t) \partial_p P^\eta_t \partial_e E^{\eta,k}_t \partial_e \cV^{\eta,k}(t,P^\eta_t,E^{\eta,k}_t) \ud t}
\end{align}
where $(\partial_p P^\eta, \partial_e E^{\eta,k})$ are the tangent processes associated to $( P^\eta,E^{\eta,k})$. They are given by
\begin{align}
\partial_e E^{\eta,k}_t &= 1 + \int_{t_0}^t \partial_y \mu(P^\eta_s,\cV^{\eta,k}(s,P^\eta_s,E^{\eta,k}_s)) \partial_e \cV^{\eta,k}(s,P^\eta_s,E^{\eta,k}_s) \partial_e E^{\eta,k}_s \ud s\,,
\\
( \partial_{p^{\ell'}}P^\eta_t)^\ell 
&= 
\1_{\set{\ell=\ell'}} + \int_{t_0}^t \partial_p b^\ell(P^\eta_s ) \partial_{p^{\ell'}}P^\eta_s \ud s \;\text{ for }\;\ell, \ell' \in \set{1,\dots,d}.
\end{align}
It is well known in our Lipschitz setting that
\begin{align}\label{eq control tangent process one more time}
\esp{\sup_{t \in [t_0,T]} |\partial_{p}P^\eta_t|^\kappa}  \le C_\kappa\,,
\end{align}
for $\kappa \ge 1$.
\\
For a positive real number $R$, we obtain from \eqref{eq starting point} that
\begin{align}
\int_{-R}^R |\partial_p \cV^{\eta,k} (t_0,p,e)| \ud e \le 
\esp{\int_{t_0}^T\int_{-R}^R|\partial_p\mu(P^\eta_t,Y^{\eta,k}_t) \partial_p P^\eta_t|   \partial_e [\cV^{\eta,k}(t,P^\eta_t,E^{\eta,k}_t)] \ud e \ud t}
\end{align}
recalling that $\partial_e E^{\eta,k}$ and $\partial_e \cV^{\eta,k}(\cdot)$ are non-negative.  Moreover, since $\cV^{\eta,k}$  is bounded and $\mu$ Lipschitz continuous, we compute
\begin{align*}
\int_{-R}^R |\partial_p \cV^{\eta,k} (t_0,p,e)| \ud e \le 2 \|\cV^{\eta,k}\|_\infty \|\partial_p\mu\|_\infty (T-t_0) \esp{\sup_{t\in[t_0,T]}| \partial_p P^\eta_t| }.
\end{align*}
Using \eqref{eq control tangent process one more time}, we eventually obtain
\begin{align}\label{eq first key bound}
\int_{-R}^R |\partial_p \cV^{\eta,k} (t_0,p,e)| \ud e \le C(T-t_0)\,.
\end{align}
Since $\cV^{\eta,k} (t_0,p+\delta,e)-\cV^{\eta,k} (t_0,p,e) = \delta \int_0^1\partial_p \cV^{\eta,k} (t_0,p+\lambda \delta,e) \ud \lambda$, we have
\begin{align}
\int_{-R}^R |\cV^{\eta,k} (t_0,p+\delta,e)-\cV^{\eta,k} (t_0,p,e)| \ud e &\le \delta \int_0^1\int_{-R}^R |\partial_p \cV^{\eta,k} (t_0,p+\lambda \delta ,e)| \ud e \lambda\,,
\\
&\le C(T-t_0) \delta
\end{align}
where we used \eqref{eq first key bound} for the last inequality. The dominated convergence theorem leads to
\begin{align}
\int_{-R}^R \frac1{\delta}|\cV(t_0,p+\delta,e)-\cV(t_0,p,e)| \ud e \le C(T-t_0) \,,
\end{align}
recalling Corollary \ref{co smooth approx of v}. Taking limit as $\delta \rightarrow 0$, we
obtain 
\begin{align}
\int_{-R}^R |\partial_p \cV(t_0,p,e)| \ud e \le C(T-t_0) \,.
\end{align}
The estimate \eqref{eq control deriv L1} is then obtained by letting $R \rightarrow +\infty$ and invoking monotone convergence.
\eproof
%\end{bluetext}

\vspace{2mm}
\noindent The study of the regularity of the decoupling field in full generality is outside the scope of this paper. We, however, give now an example of model for which Assumption \ref{ass restrict regul} is known to hold.
{\color{black}
\begin{Remark}\label{re example for ass} Let us consider the following \emph{toy model}, which has been studied in \cite{carmona2013singularfirst}
\begin{equation}\label{eq toy model}
\left\{
\begin{split}
dP_t& = \bar{\sigma}\ud \cW_{t}\,, \; P_0 = p\,, \\%\label{eq toy model P}\\
dE_t &= \left(P_t-Y_t\right)\ud t\,, \; E_0 = e\,, \\%\label{eq toy model E}\\
dY_t &=  Z_t  \ud \cW_t. \\%\label{eq toy model Y}
\end{split}
\right.
\end{equation}
Here, $\bar{\sigma}$ is a positive  constant, $\cW$ is a one-dimensional Brownian Motion and the terminal condition (in the weak sense of \eqref{relaxed terminal condition with p}) for $Y$ is given by $\phi(E_T)$, recall \eqref{eq very special spec phi}. By Theorem \ref{th existence uniqueness one-period}  above, it has a unique solution and we denote by $\upsilon$ its decoupling field. In particular, we have $\upsilon(0,p,e) = Y_0$. Now, we introduce, as in \cite{carmona2013singularfirst}, the following process
\begin{align}
\bar{E}_t = E_t + (T-t)P_t \,,
\end{align}
and  observe then that $(\bar{E},Y)$ satisfies
\begin{equation}\label{eq toy model reduced}
\left\{
\begin{split}
\ud \bar{E}_t &= - Y_t \ud t + (T-t) \ud \cW_t
\\
\ud Y_t &= Z_t \ud \cW_t
\end{split}
\right.
\end{equation}
with the same terminal condition since $\bar{E}_T = E_T$.
Independently from \eqref{eq toy model}, the system \eqref{eq toy model reduced} is thoroughly studied in \cite[Section 4]{carmona2013singularfirst}, this system has also a unique solution with a decoupling field $Y_0 = \bar{\upsilon}(0,\bar{E}_0)$, and from Proposition 6 in \cite{carmona2013singularfirst} again, $\bar{\upsilon}$ is $\mathscr{C}^{1,2}([0,T)\times \R)$. We deduce that
$\upsilon(0,p,e)=\bar{\upsilon}(0,e+Tp)$ and more generally $\upsilon(t,p,e)=\bar{\upsilon}(t,e+(T-t)p)$ for $t\in [0,T)$. This yields that $\upsilon$ is $\mathscr{C}^{1,2,2}([0,T)\times\R^d\times\R)$. 
%Moreover, one computes that $\partial_p(t,p,e)\upsilon= (T-t)\partial_\xi \bar{\upsilon}(t,e+(T-t)p)$ which shows, from???, that $\partial_p \upsilon$ is bounded.
\end{Remark}
}%ec

\vspace{2mm}
\noindent Let us observe that Proposition \ref{pr property of V smooth} states estimates for the first order derivatives of $\cV$. And, though $\cV$ is assumed to be smooth in Assumption \ref{ass restrict regul}, nothing is said on the control of higher order derivatives. As we will see later on, the control of  the second order derivative and third order derivative with respect to $p$, will be needed. This should not come as a surprise as we will approximate the $P$ process by a discrete in time and space process, see Section \ref{subsubse approx generic}. We shall not make further assumption on the behavior of $\cV$ and its derivatives. Instead, we introduce a proxy to $\cV$ that will be used in the proof of convergence of our numerical algorithm.

\vspace{2mm}
\noindent Let then $\cV^{\epsilon}$ be a smoothing  of the decoupling field $\cV$, namely: For $\epsilon \in (0,1)$,
\begin{align}\label{eq de cV epsilon}
	\cV^\epsilon(t,p,e) := \int \cV(t,p+q,e) \varphi_\epsilon(q) \ud q \quad\text{ where } \; \varphi_\epsilon(q) := \frac1\epsilon\varphi(\frac{q}\epsilon)\,,\, q \in \R,
	\end{align}
	and with $\varphi(\cdot)$ a smooth compactly supported probability density function.
\vspace{2mm}

\noindent We now show that $\cV^\epsilon$ satisfies the same PDE as $\cV$ up to an error term that we can control. From the previous Assumption \ref{ass restrict regul}, we first observe that $\cV^\epsilon$, defined in
\eqref{eq de cV epsilon} belongs also to $\mathscr{C}^{1,2,1}([0,T) \times \R^d \times \R)$ and that it satisfies the same bounds as $\cV$, namely
\begin{align}\label{eq control deriv epsilon}
0 &\le \partial_e \cV^\epsilon(t,p,e) \le \frac{C}{T-t}\;, \quad |\partial_p \cV^\epsilon(t,p,e)|  \le C
\\
&  \text{ and }
 \int |\partial_p \cV^\epsilon(t,p,z)| \ud z \le C(T-t)\;, \quad \forall (t,p,e) \in [0,T)\times \R^d \times \R\;.
\end{align}

\begin{Proposition}
Under Assumption \ref{ass restrict regul}, the function $\cV^{\epsilon}$ satisfies on $[0,T)\times\R^d\times \R$,
\begin{align} \label{eq pde regul regul V}
\partial_t \cV^{\epsilon} &+b(p)\partial_p \cV^{\epsilon}+ \frac12\bar{\sigma}^2\Delta_p \cV^{\epsilon}
%\\
%&
+  \mu(\cV(t,p,e), p) \partial_e \cV^{\epsilon}  = \theta^\epsilon(t,p,e)
\end{align}
with
\begin{align}\label{eq control theta}
|\theta^\epsilon(t,p,e)|
\le C\epsilon\left (\int |\partial_p \cV^{\epsilon}(t,p+q,e)|\varphi_\epsilon(q) \ud q  + \partial_e\cV^{\epsilon}(t,p,e)\right )\,.
\end{align}
%\textcolor{red}{simplifier encore l'expression du controle ci-dessus}
%We also observe that $\cV^\epsilon$ satisfies the same bounds as $\cV$ given in \eqref{eq control deriv}.
%\textcolor{red}{We will need to control
%\begin{align}
%\int \int |\partial_p \cV^{\eta,k}(t,p+q,e)|\varphi_\epsilon(q) \ud q \ud e
%\end{align}
%this is due to the non linearity of $b$.
%}
\end{Proposition}

%\begin{bluetext}
\proof
Recall that $\cV$ satisfies \eqref{eq smooth version},
%\begin{align}
%\partial_t \cV^{\eta,k,\epsilon} +{\mu(\cV^{\eta,k,\epsilon},p)} \partial_{e} \cV^{\eta,k,\epsilon} + b(p)\partial_p \cV^{\eta,k,\epsilon} + \mathrm{L}^{\!\eta}\cV^{\eta,k,\epsilon} = 0
%\end{align}
and that, for $(t,p,e)\in [0,T)\times\R^d\times \R$,
\begin{align}
\partial_t \cV^{\epsilon}(t,p,e) &= \int \partial_t \cV(t,p+q,e) \varphi_\epsilon(q) \ud q \text{ and }
\partial^2_{pp} \cV^{\epsilon}(t,p,e) = \int \partial^2_{pp}\cV(t,p+q,e) \varphi_\epsilon(q) \ud q.
\end{align}
We thus compute, by linearity of $\partial_t$, $\Delta_{p}$,
\begin{align}
\partial_t \cV^{\epsilon}(t,p,e) &+ \frac12 \bar{\sigma}^2 \Delta_{p} \cV^{\epsilon}(t,p,e)
= -\int b(p+q)\partial_p \cV(t,p+q,e)  \varphi_\epsilon(q) \ud q
\\
&
- \int \mu(\cV(t,p+q,e), p+q) \partial_e \cV(t,p+q,e) \varphi_\epsilon(q) \ud q
%
%\\
%{F(u,p)} \partial_{e} u + \cL_p u + \frac12\eta^2 (\partial^2_{ee}u+ \Delta_{pp}u) 
\end{align}
which leads to
\begin{align}
\partial_t \cV^{\epsilon}(t,p,e) &+b(p)\partial_p \cV^{\epsilon}(t,p,e)+  \frac12\bar{\sigma}^2 \Delta_{p} \cV^{\epsilon}(t,p,e)
\\
&+  \mu(\cV(t,p,e), p) \partial_e \cV^{\epsilon}(t,p,e)  = \theta^\epsilon(t,p,e)
\end{align}
with
\begin{align}
\theta^\epsilon(t,p,e) &:= \int \set{b(p) - b(p+q)}\partial_p \cV(t,p+q,e) \varphi_\epsilon(q)  \ud q
\\
&+ \int \set{\mu(\cV(t,p,e), p)-\mu(\cV(t,p+q,e), p+q)}\partial_e \cV(t,p+q,e) \varphi_\epsilon(q) \ud q.
\end{align}
We compute, using the uniform Lipschitz property of $b$, $\mu$ and $\cV$ in the $p$ variable and the non-negativity of $\partial_e\cV$, that
\begin{align}
|\theta^\epsilon(t,p,e)| &\le C \int|q||\partial_p \cV(t,p+q,e)|\varphi_\epsilon(q) \ud q + C\int |q| \partial_e\cV(t,p+q,e) \varphi_\epsilon(q) \ud q
\\
& \le C\epsilon\left (\int |\partial_p \cV(t,p+q,e)|\varphi_\epsilon(q) \ud q  + \partial_e\cV^{\epsilon}(t,p,e)\right )
\end{align}
where we use the fact that $\varphi$ has compact support to get the last inequality.
\eproof
%\end{bluetext}

To conclude this section, we collect some useful control in the $L^1$-sense of the function $\cV^\epsilon$ and its derivatives. Their proof is very classical and we skip it.

\begin{Lemma}\label{le bound for V epsilon}
Assumption \ref{ass restrict regul}, the followings hold, for $(t,p) \in [0,T)\times\R^d$,
\begin{align}\label{eq control  dppp u epsilon}
  \epsilon\int |\partial^2_{p_ip_j} \cV^\epsilon(t,p,e)| \ud e + \epsilon^2\int |\partial^3_{p_ip_jp_k} \cV^\epsilon(t,p,e)| \ud e   \le C\;,
\end{align}
for $i,j,k \in \set{1,\dots,d}$.
And thus, for $(t,p,p') \in [0,T)\times\R^d\times\R^d$,
\begin{align}\label{eq control lip dpp u epsilon}
 \int |\partial_{p} \cV^{\epsilon}(t,p,e) -\partial_{p} \cV^{\epsilon}(t,p',e) |\ud e + \epsilon \int |\partial^2_{pp} \cV^{\epsilon}(t,p,e) -\partial^2_{pp} \cV^{\epsilon}(t,p',e) |\ud e \le \frac{C}{\epsilon}|p-p'|,
\end{align}
where $C$ is a positive constant.
\end{Lemma}

\section{Numerical algorithm}
\label{se algo}

In this section, we first introduce various versions of the splitting scheme that will be implemented in practice.
We  then introduce the  numerical errors that have to be taken into account in order to obtain the convergence of the scheme. The main error decomposition is given at the end of this section, see Proposition \ref{pr global stab}, together with our main convergence result, see Theorem \ref{th main conv result}. Precise study of each error is postponed to the next section.

\subsection{ Fully implementable schemes}
\label{subse full scheme}

The schemes we consider are all based on the approximation of the Wiener measure by a  discrete and finite set of path. Indeed, on the grid $\pi$, the Brownian motion is approximated by a discrete time process denoted $\widehat{W}$. At each date $t_n$, ${\mathfrak{S}}_n \subset \R^d$  denotes the support of the law of the random variable $\widehat{W}_{t_n}$. To define $\widehat{W}$, we first define $(\Delta \widehat{W}_n := \widehat{W}_{t_{n+1}} - \widehat{W}_{t_n})_{0 \le n \le N-1}$ which stands for discrete approximation of the  Brownian increments $(W_{t_{n+1}} - W_{t_n})_{0 \le n \le N-1}$.

\begin{Definition}[Brownian increments approximation] \label{de bro inc}
 We use  the cubature formula introduced in \cite[Section A.2]{gyurko2011efficient} which requires only $2d$ points. Denoting $(\mathfrak{e}^\ell)_{1\le \ell \le d}$ the canonical basis of $\R^d$, we set, for $1 \le i \le I=2d$,
$\P(\Delta \widehat{W}_n=\omega^i_{\mathfrak{h}}) = \frac1{2d}$ and $\omega^i_{\mathfrak{h}} = -\sqrt{d{\mathfrak{h}}}\mathfrak{e}^\ell$, if $i=2\ell$ or $\omega^i_{\mathfrak{h}} = \sqrt{d{\mathfrak{h}}}\mathfrak{e}^\ell$ if $i=2\ell - 1$. 
\end{Definition}

\noindent By construction, we observe that 
\begin{align}\label{eq properties of widehat W}
 \text{ for } x \in \mathfrak{S}_n, \, x + \Delta \widehat{W}_n \in \mathfrak{S}_{n+1}\,.
\end{align}
Moreover, the set $\mathfrak{S}_n$ are discrete and finite.

\subsubsection{Schemes for generic diffusion}
\label{subsubse approx generic}
In this section, we treat the general case where $P$ is solution to a SDE with Lipschitz coefficients as in \eqref{eq singular fbsde}.
The discrete approximation on the grid $\pi$ is given by the classical Euler Scheme
\begin{align}\label{de disc euler scheme P}
\widehat{P}_0 = P_0 \text{ and } \widehat{P}_{t_{n+1}} = \widehat{P}_{t_{n}} + b(\widehat{P}_{t_{n}}) \mathfrak{h} + \sigma(\widehat{P}_{t_{n}}) \Delta \widehat{W}_n\,,
\end{align}
where $(\Delta \widehat{W}_n)$ is introduced in Definition \ref{de bro inc}.
For later use, we also define, for $p \in \R^d$ and any $0 \le n \le N-1 $,
\begin{align}\label{eq de one step disc euler scheme P}
\widehat{P}^{t_n,p}_{t_{n+1}} := p +  b(p) \mathfrak{h} + \sigma(p) \Delta \widehat{W}_n\,,
\end{align}
and, for $1\le i \le 2d $,
\begin{align}\label{eq de one step disc euler scheme P one value}
\left(\widehat{P}^{t_n,p}_{t_{n+1}} \right)^i := p +  b(p) \mathfrak{h} + \sigma(p) \omega^i_{\mathfrak{h}}\,.
\end{align}
Let us denote by $\cP_n$ the discrete and finite support of $\widehat{P}_{t_{n}}$  for $0 \le n \le N-1$. We observe, due to the very definition of 
$(\widehat{P}_{t_n})$, that
\begin{align}
\text{ for } p \in \cP_n,\; \widehat{P}^{t_n,p}_{t_{n+1}} \in \cP_{n+1}\;.
\end{align}

\noindent In this context, we first  introduce a discrete version of the operator $\cT$, recall Definition \ref{de transport step}, that will compute an approximation to \eqref{eq de cons law associated with transport} written in \emph{forward form}: We shall use the celebrated Sticky Particle Dynamics (SPD) \cite{brenier1998sticky}  see also \cite[Section 1.1]{jourdainreygner16}. The SPD is particularly simple to implement in our case, since, due to the monotonicity assumption on $(\mu,\psi)$, there is no particle colliding!  

\noindent For $M \ge 1$, let 
\begin{align}
D_M &:=\set{\mathrm{e}=(e_1,\dots,e_m,\dots,e_M) \in \R^M \,|\, e_1 \le \dots \le e_m \le \dots \le e_M}
\\
\text{and} \;\mathcal{I}_M &:= \set{\theta \in \mathscr{I} | \theta(\cdot) := {H}\!*\!(\frac1M \sum_{m=1}^M \delta_{e_m}),\mathrm{e}\in D_M}
\end{align}
 where $H$ is the Heaviside function $x \mapsto \1_{\set{x \ge 0}}$, $*$ the convolution operator and $\delta_e$ is the Dirac mass at $e \in \R$.
 \\
 For later use, we introduce, for $M \ge 0$ and $n \in \set{0,\dots,N}$, the set of function
 \begin{align}\label{eq de KnM}
 \cK_{n}^M  = \set{\psi:\cP_n \times \R \rightarrow [0,1]\,|\, \text{ for } p \in \cP_n,\, \psi(p,\cdot) \in \cI_M}\;.
 \end{align}
 \\
 The discrete version of $\cT$, denoted $\mathfrak{T}$,  acts on empirical CDF belonging to $\mathscr{I}_M$
 and is given by
\begin{align} \label{de disc op T}
(0,\infty) \times \R^d \times \mathcal{I}_M \ni (h,p,\theta) \mapsto \mathfrak{T}_h^{M}(p,\theta):= {H}\!*\!(\frac1M \sum_{m=1}^M \delta_{e^p_m(h)}) \in \mathcal{I}_M\,.
%\\
%=(E^{p,m}_h)_{1 \le m \le M} \in D_M\;,
\end{align}

\noindent Above $(e^{p}_m(h))_{1 \le m \le M}$ is a set (of positions) of particles computed as follows.
Given the initial position $\mathrm{e}(0) \in D_M$ (representing $\theta$) and velocities $(\bar{F}_m(p))_{1 \le m \le M}$ set to 
$\bar{F}_m(p) = - \int_{(m-1)/M}^{m/M}\mu(p,y)\ud y$, we consider $M$ particles $(e^{p}_m)_{1 \le m \le M}$, whose positions at time $t \in [0,h]$ are simply given by
\begin{align}\label{eq particle dynamics}
e^{p}_m(t) = e_m(0) + \bar{F}_m(p) t\;.
\end{align}
We observe that $(e^{p}_m(t))_{1 \le m \le M} \in D_M$, for all $t \in [0,h]$, as $-\mu$ is non-decreasing.

%We will still rely on the discrete operator $\mathfrak{T}^M_{\mathfrak{h}}$ given in \eqref{de disc op T}. 
\vspace{2mm}
We now present the approximation used for the approximation of $\cD_{\mathfrak{h}}(\psi)$, recall Definition \ref{de diffusion step} and \eqref{eq other de bar v}.
 
\noindent Let $\psi \in \cK^M_{n+1}$ and fix $\mathrm{p} \in \cP_n$. We observe 
\begin{align}
\bar{v}_n(\mathrm{p},e) &:=\esp{\psi(\widehat{P}^{t_n,\mathrm{p}}_{t_{n+1}},e) } 
%\nonumber
%\\
%&
= \frac1{2d} \sum_{i=1}^{2d} \psi(\left(\widehat{P}^{t_n,p}_{t_{n+1}} \right)^i,e) \label{eq clarify diff generic}
\end{align}
recall \eqref{eq de one step disc euler scheme P one value}.
\\
Since $\psi \in \cK^M_{n+1}$, recall \eqref{eq de KnM}, we know that there is some $\mathrm{e}^i \in D_M$, such that
\begin{align*}
\psi(\left(\widehat{P}^{t_n,p}_{t_{n+1}} \right)^i,e) = H\!*\!(\frac1M \sum_{m=1}^M \delta_{e^{i}_m})(e)\,.
\end{align*}
By linearity, \eqref{eq clarify diff generic} reads
\begin{align}\label{eq def bar v}
\bar{v}_n(\mathrm{p},e) = H\!*\!(\frac1{2dM} \sum_{i=1}^I \sum_{m=1}^M \delta_{e^i_m})(e)\,.
\end{align}
The approximation of the diffusion operator is thus given by
\begin{align}\label{eq def approx op D}
\cK_{n+1}^M \ni \psi \mapsto \mathfrak{D}_n^M(\psi) :=\bar{v}_n  \in \cK_{n}^{2Md}\,.
\end{align}
Finally, the scheme will have essentially two versions:
\begin{itemize}%\label{two diff cases}
	\item[\textbf{CASE 1}:] We keep the $2dM$ particles at each step $n$. The overall procedure will then be the iteration of the two operators $\mathfrak{T}$ and $\mathfrak{D}$, see Definition \ref{de art scheme case 1} below.
%	The approximation operator $\mathfrak{D}_n^M$ is defined by
%	\begin{align}\label{eq de particle regression definition1}
%	\cK_{n+1}^M \ni \psi \mapsto \mathfrak{D}_n^M(\psi) :=\bar{v}  \in \cK_{n}^{2Md}\;.
%	\end{align}
\item[\textbf{CASE 2}:] There is no need to keep $2dM$ particles at step $n$, when the function $\psi$ at step $n+1$ is given by $M$ particles (for each $p \in \cP_{n+1}$). 
To reduce the number of particles, we apply another operator $\mathfrak{R}$, namely, for $M \ge m\ge 1$,
\begin{align}\label{eq formal de R}
\cI^{M} \ni \psi \mapsto \mathfrak{R}^{M,m}(\psi)   \in \cI^m \;.
\end{align}
Various implementation are possible, we refer to the numerical section for a precise description. The overall procedure for this case is given in Definition \ref{de art scheme case 2} below.
%first sort the cloud of particles $\cE$ to obtain $\tilde{\mathrm{e}}^{\mathrm{w}} \in D_{2dM}$, then we consider $\bar{\mathrm{e}}^{\mathrm{w}}:=(\tilde{e}^{\mathrm{w}}_{2dm})_{1 \le m \le M}$.
%The approximation operator $\mathfrak{D}_n^M$ is finally defined by %\textcolor{magenta}{}
%\begin{align}\label{eq de particle regression definition}
%\cK_{n+1}^M \ni \psi \mapsto \mathfrak{D}_n^M(\psi) = \textcolor{red}{\tilde{v}} \in 
%\cK_{n}^M\;.
%\end{align}

\end{itemize}

\vspace{2mm} 
\noindent Let us now finally introduce the scheme formally.

\begin{Definition}[Scheme in CASE 1] \label{de art scheme case 1} Fix $M \ge 1$ and set for $0 \le n \le N$, $M_n := (2d)^{N-n}M$.
\begin{enumerate}
\item  At $n=N$: Set $\mathrm{e}_N :=(\Lambda,\dots,\Lambda)\in D_{M_N}$ whose empirical CDF is the terminal condition $\phi = \1_{\set{e \ge \Lambda}}$, recall \eqref{eq very special spec phi}. Setting simply ${u}^{N,M}_N(p,\cdot) = \phi$, for all $p \in \cP_{N}$, we observe that ${u}^{N,M}_N \in \cK_N^{M_N}$.

\item For $n<N$: Given ${u}^{N,M}_{n+1} \in \cK_{n+1}^{M_{n+1}}$, define
 \begin{align}\label{eq de bar u NM}
 \bar{u}^{N,M}_n =  \mathfrak{D}_n^{M_{n+1}}({u}^{N,M}_{n+1}) \in \cK_{n}^{M_{n}}\,,
 \end{align}
 recall \eqref{eq def approx op D}, and
then ${u}^{N,M}_n$ by, for each $p \in \cP_n$,
\begin{align}
{u}^{N,M}_n(p,\cdot) = \mathfrak{T}_{\mathfrak{h}}^{M_{n}}(p,\bar{u}^{N,M}_n(p,\cdot)) \,,
\end{align}
recall \eqref{de disc op T}.
\end{enumerate}
The approximation of $\cV(0,P_0,\cdot)$ is then given by $u^{N,M}_0(P_0,\cdot)$. 
 \end{Definition}
 
\begin{Definition}[Scheme in CASE 2] \label{de art scheme case 2} Fix $M \ge 1$.
\begin{enumerate}
\item  At $n=N$: Set $\mathrm{e}_N :=(\Lambda,\cdots,\Lambda)\in D_M$ whose empirical CDF is the terminal condition $\phi = \1_{\set{e \ge \Lambda}}$, recall \eqref{eq very special spec phi}. Setting simply $v^{N,M}_N(p,\cdot) := \phi$, for all $p \in \cP_{N}$, we observe that $v^{N,M}_N \in \cK_N^M$.

\item For $n<N$: Given $v^{N,M}_{n+1} \in \cK_{n+1}^M$, define
 $\bar{v}^{N,M}_n =  \mathfrak{D}_n^M(v^{N,M}_{n+1}) \in \cK_{n}^{2dM}$,  recall \eqref{eq def approx op D}, and
then for each $p \in \cP_n$, 
\begin{align}
\check{v}^{N,M}_n(p,\cdot) =  \mathfrak{R}^{2dM,M}(\bar{v}^{N,M}_n(p,\cdot))
\end{align}
recall \eqref{eq formal de R}, and finally ${v}^{N,M}_n(p,\cdot) $ by, 
\begin{align}
v^{N,M}_n(p,\cdot) = \mathfrak{T}_{\mathfrak{h}}^M(p,\check{v}^{N,M}_n(p,\cdot))\,,
\end{align}
recall \eqref{de disc op T}.
\end{enumerate}
The approximation of $\cV(0,P_0,\cdot)$ is then given by $v^{N,M}_0(P_0,\cdot)$. 
 \end{Definition}

\subsubsection{A simplified functional Brownian setting}
\label{subsubse brownian setting}
We now consider a special case for the process $P$ for which the numerical implementation is less computationally demanding, see Remark \ref{re blabla tree recombining} below.
\begin{Assumption} \label{ass brownian functional}
The process $P$ is given as a function of the Brownian motion, namely
\begin{align*}
P_t = \mathfrak{P}(t,W_t) \;\text{ for }\; \mathfrak{P}:[0,T]\times\R^d \mapsto \R^d\;.
\end{align*}
\end{Assumption}

\begin{Remark} The main application of this case is when $P$ is the Brownian motion itself (also known as the Bachelier model in financial application) or a Geometric Brownian Motion (also known as Black-Scholes model in financial application).
\end{Remark}
%
%On the grid $\pi$, the Brownian motion is approximated by a discrete time process denoted $\widehat{W}$ . At each date $t_n$, ${\mathfrak{S}}_n \subset \R^d$  denotes the support of the law of the random variable. We assume that it is discrete and finite. In practice, to define $\widehat{W}$, we first define $(\Delta \widehat{W}_n := \widehat{W}_{t_{n+1}} - \widehat{W}_{t_n})_{0 \le n \le N-1}$ which stands for discrete approximation of the  Brownian increments $(W_{t_{n+1}} - W_{t_n})_{0 \le n \le N-1}$. 
%
%
%\begin{Definition}[Brownian increments] \label{de bro inc}
% We use  the cubature formula introduced in \cite[Section A.2]{gyurko2011efficient} which requires only $2d$ points. Denoting $(\mathfrak{e}^\ell)_{1\le \ell \le d}$ the canonical basis of $\R^d$, we set, for $1 \le i \le I=2d$,
%$\P(\Delta \widehat{W}_n=\omega^i_{\mathfrak{h}}) = \frac1{2d}$ and $\omega^i_{\mathfrak{h}} = -\sqrt{d{\mathfrak{h}}}\mathfrak{e}^\ell$, if $i=2\ell$ or $\omega^i_{\mathfrak{h}} = \sqrt{d{\mathfrak{h}}}\mathfrak{e}^\ell$ if $i=2\ell - 1$. 
%\end{Definition}
%
%
%We observe that 
%\begin{align}\label{eq properties of widehat W}
% \text{ for } x \in \mathfrak{S}_n, \, x + \Delta \widehat{W}_n \in \mathfrak{S}_{n+1}\,.
%\end{align}
\noindent Compared to the previous section, the only difference is the approximation of $P$. In this context,  it is naturally given by
\begin{align}\label{eq de hat P}
\widehat{P}_{t} := \mathfrak{P}(t,\widehat{W}_{t})\,,\; t \in \pi\,.
\end{align}
Note that we slightly abuse the notation here by keeping the same ones as the ones introduced in the previous section.

\noindent So, for $0 \le n \le N$, 
\begin{align}
\cP_{n}=\set{\mathrm{p} = \mathfrak{P}(t_{n},\mathrm{w}),\mathrm{w} \in \mathfrak{S}_{n}}\,,
\end{align}
which is  the (discrete) support of $\widehat{P}_{t_{n}}$.
We also define, for $p =  \mathfrak{P}(t_{n},\mathrm{w})$,
\begin{align}
\widehat{P}^{t_n,p}_{t_{n+1}} = \mathfrak{P}(t_{n+1},\mathrm{w}+\Delta \widehat{W}_n) \text{ and }
\left(\widehat{P}^{t_n,p}_{t_{n+1}}\right)^i= \mathfrak{P}(t_{n+1},\mathrm{w}+ \omega^i_{\mathfrak{h}})\,,
\end{align}
for $1 \le i \le 2d$, recall Definition \ref{de bro inc}.

\vspace{2mm}
\noindent Now that the approximation of $P$ has been clarified in this context, we can use the various schemes given in Definition \ref{de art scheme case 1} and Definition \ref{de art scheme case 2}. 
\begin{Remark}\label{re blabla tree recombining}
The main difference between the generic diffusion case and the functional Brownian setting comes from the numerical complexity associated to their implementation. Indeed, it is the case that the growth size of $\mathfrak{S}_n$ with respect to $n$ is tamed (for low dimensional problem) because the ``tree'' associated to the brownian motion is naturally recombining. It is a priori not the case for generic diffusion and in particular the growth of $\cP_n$ is exponential in $n$. For implementation purposes, we just use the Brownian setting. A  way to control the growth for generic diffusion is to introduce some interpolation on a grid as e.g. in \cite{chassagneux2020cubature}, this method is well known and we do not present it here. It introduces a further error that would need to be controlled. Note that all the methods  we have described here and we study, are impacted by the ``curse of dimensionality'' associated to the dimension of $P$. However, we should also note that for dimension around 5 and in the case of the brownian setting, the problems are tractable numerically, see section \ref{se numerics}.
\end{Remark}

\subsection{Error decomposition and main result}

%For ease of presentation, we define 
%\begin{align}
%u^{N,M}(t_n,p,e) = \gamma_n(p,e) \;\text{ for }\; p \in \cP_n\;,
%\end{align}
%where $\gamma$ is the solution of the scheme defined above.
%\\

We first present the various sources of errors introduced by the schemes above. In the sequel, we shall mainly conduct our analysis for the scheme given in \textsc{CASE 1}.

\vspace{3mm}

\noindent The error we seek to control is, for $\gamma_n = u^{N,M}_n$ or $\gamma_n = v^{N,M}_n$, recall Definition \ref{de art scheme case 1} and \ref{de art scheme case 2},
\begin{align}\label{eq main error}
\mathrm{err}(N,M) := \int |\gamma_0(P_0,e)-\cV(0,P_0,e)|\ud e\;.
\end{align}	
	
\noindent And we first observe
	\begin{align}
	\mathrm{err} \le {\cE}_0(P_0) + {\cE}^r_0\,,
	\end{align}
	with for $0 \le n \le N$ and $p \in \cP_n$,
	\begin{align}
	{\cE}^r_n (p)&:= \int |\cV^{\epsilon}(t_n,p,e)-\cV(t_n,p,e)|\ud e  \,,
	\\
%	\end{align}
%	and 
%\begin{align}
{\cE}_n(p)  &:= \int |\gamma_n(p,e)-\cV^\epsilon(t_n,p,e)|\ud e \,.
\end{align}

\noindent In \textsc{CASE 1}, we  define, for ${p} \in \cP_n$,
	\begin{align}
	&\cE^{\mathrm{T}}_n(p) := \int \Big|\mathfrak{T}^{M_n}_{\mathfrak{h}}(p,\bar{u}^{N,M}_n (p,\cdot))(e) - \tilde{\cT}_h(p,\bar{u}^{N,M}_n(p,\cdot))(e)\Big|\ud e\,, \label{eq error def T} \\
	&\bar{\cE}_n(p):= \int \Big|\tilde{\cT}_h(p,\bar{u}^{N,M}_n(p,\cdot))(e) - \tilde{\cT}_h(p,\mathbb{E}[\cV^{\epsilon}(t_{n+1},\widehat{P}_{t_{n+1}}^{t_{n},p},\cdot)])(e)\Big|\ud e \,, \label{eq error def stab}
	\end{align}
recall Definition \ref{de art scheme case 1} and Remark \ref{re abuse transport op}.\\
	Moreover, for $p \in \R^d$, we set
	\begin{align}
	& \cE^{D}_n(p):= \int \Big|\tilde{\cT}_h(p,\mathbb{E}[\cV^{\epsilon}(t_{n+1},\widehat{P}_{t_{n+1}}^{t_{n},p},\cdot)])(e) - \tilde{\cT}_h(p,\mathbb{E}[\cV^{\epsilon}(t_{n+1},{P}_{t_{n+1}}^{t_{n},p},\cdot)])(e) \Big|\ud e\,, \label{eq error def diff} \\ 
	& \cE^{S}_n(p):= \int \Big|\tilde{\cT}_h(p,\mathbb{E}[\cV^{\epsilon}(t_{n+1},{P}_{t_{n+1}}^{t_{n},p},\cdot)])(e)- \cV^{\epsilon}(t_{n},p,e)\Big|\ud e\,. \label{eq error def split}
	\end{align}
We have the following.
\begin{Lemma} \label{le one step decomp}
Under \textsc{CASE 1}, for $0 \le n \le N-2$,
\begin{align}
{\cE}_n(p) \le \cE^{\mathrm{T}}_n(p) + \bar{\cE}_n(p) + \cE^{D}_n(p) + \cE^{S}_n(p)
\end{align}
and
\begin{align}\label{eq control last step}
{\cE}_{N-1}(p)  &\le \cE^{\mathrm{T}}_{N-1}(p)   + \cE^{\mathrm{r}}_{N-1}(p) 
+
\int | \cS_{\mathfrak{h}}(\phi(\cdot))(p,e) - \Theta_{\mathfrak{h}}(\phi(\cdot))(p,e)| \ud e %\nonumber
%\\
%&\quad +
%\int | \cV(t_{N-1},p,e) - \cV^\epsilon(t_{N-1},p,e)| \ud e \;.
\end{align}
observing that
\begin{align} \label{eq last step transport approx}
\cE^{\mathrm{T}}_{N-1}(p) = \int |\mathfrak{T}^{M}_{\mathfrak{h}}(p,\phi(\cdot))(e) - \cT_{\mathfrak{h}}(p,\phi(\cdot)(e))| \ud e\,.
\end{align}
\end{Lemma}

%\begin{bluetext}
\proof
The first inequality is straightforward. For the second, we observe that $u^{M,N}_N = \phi$ and does not depend on $p$ variable so that $\bar{u}^{M,N}_{N-1} = \phi$. %We similarly observe that $\cV^\epsilon(T,\cdot)=\phi(\cdot)$.
\begin{align*}
{\cE}_{N-1}(p)  &:= \int |u^{M,N}_{N-1}(p,e)-\cV^\epsilon(t_{N-1},p,e)|\ud e
\\
&= \int |\mathfrak{T}^{M}_{\mathfrak{h}}(p,\phi(\cdot))(e)-\cV^\epsilon(t_{N-1},p,e)|\ud e
\\
&\le  \int |\mathfrak{T}^{M}_{\mathfrak{h}}(p,\phi(\cdot))(e)- \cT_{\mathfrak{h}}(p,\phi(\cdot))(e)| \ud e
+
\int | \cT_{\mathfrak{h}}(p,\phi(\cdot))(e) - \cV(t_{N-1},p,e)| \ud e
\\
&\quad +
\int | \cV(t_{N-1},p,e) - \cV^\epsilon(t_{N-1},p,e)| \ud e.
\end{align*}
We then observe that $\cT_{\mathfrak{h}}(p,\phi(\cdot))(e) = \cS_{\mathfrak{h}}(\phi(\cdot))(p,e)$ as again $\phi$
does not depend on $p$.
\eproof
%\end{bluetext}

\vspace{2mm}
\noindent We now comment quickly the various errors appearing above:
The term $\cE^{\mathrm{T}}_n(p)$ is the local error due to the approximation of the transport operator by the SPD. The term $\bar{\cE}_n(p)$ will allow to propagate the local errors thanks to the $L^1$ stability of $\cT_h(\cdot)$. The term $\cE^{D}$ is linked to the approximation of the diffusion by the tree and $\cE^{S}$ is the splitting error applied to the proxy $\cV^\epsilon$.

\begin{Proposition}[Stability]\label{pr global stab}
Under \textsc{CASE 1}, the following hold
\begin{align}
{\cE}_0(P_0) \le \sum_{n=0}^{N-2} \esp{ \cE^{\mathrm{T}}_n(\widehat{P}_{t_n})   +  \cE^{D}_n(\widehat{P}_{t_n}) +  \cE^{S}_n(\widehat{P}_{t_n}) }  + \esp{{\cE}_{N-1}(\widehat{P}_{t_{N-1}})}\;.
\end{align}
\end{Proposition}

%\begin{bluetext}
\proof Since $\tilde{\cT}_h$ is $L^1$-stable, \eqref{eq error def stab} yields for $p \in \cP_n$, 
\begin{align*}
\bar{\cE}_n(p) \le \int |\bar{u}^{N,M}(t_n,p,e)-\mathbb{E}[\cV^{\epsilon}(t_{n+1},\widehat{P}_{t_{n+1}}^{t_{n},p},e)]| \ud e.
\end{align*}
Under \textsc{CASE 1}, recall \eqref{eq de bar u NM}, we then  get
\begin{align}
\bar{\cE}_n(p) \le \esp{\int |{u}^{N,M}(t_{n+1},\widehat{P}_{t_{n+1}}^{t_{n},p},e)- \cV^{\epsilon}(t_{n+1},\widehat{P}_{t_{n+1}}^{t_{n},p},e)| \ud e}.
\end{align}
Thus,
\begin{align}
\esp{\bar{\cE}_n(\widehat{P}_{t_n})} \le \esp{{\cE}_{n+1}(\widehat{P}_{t_{n+1}})}.
\end{align}
We thus  deduce from Lemma \ref{le one step decomp}
\begin{align}
\esp{ {\cE}_n(\widehat{P}_{t_n}) } \le \esp{ \cE^{\mathrm{T}}_n(\widehat{P}_{t_n})   +  \cE^{D}_n(\widehat{P}_{t_n}) +  \cE^{S}_n(\widehat{P}_{t_n}) } +  \esp{{\cE}_{n+1}(\widehat{P}_{t_{n+1}})}.
\end{align}
The proof is concluded by induction observing $\esp{{\cE}_{N}(\widehat{P}_{t_{N}})} = 0$.
\eproof
%\end{bluetext}

\vspace{2mm}
The next section is dedicated to the study of the various error appearing above. We prove in Section \ref{subse proof th main conv result} below the main theoretical result of this work:

\begin{Theorem}\label{th main conv result}
Let Assumption \ref{ass restrict regul} hold. Then, under \textsc{CASE 1}, 
\begin{align}\label{eq control main error}
\mathrm{err}(N,M) \le C (\frac1M + \frac{\sqrt{\mathfrak{h}}}{\epsilon^2} + \epsilon )\,,
\end{align}
recall \eqref{eq main error}.
Moreover, setting $\epsilon = \mathfrak{h}^{\frac16}$ and $M = \frac1\epsilon$, we have
\begin{align}\label{eq main err in term of h}
 \int|\cV(0,P_0,e) - u^{N,M}_0(P_0,e)| \ud e\le C \mathfrak{h}^{\frac16}.
\end{align}

\end{Theorem}

\section{Study of the errors}
\label{se errors}
\subsection{Approximation of transport operator}
%In the above, we have introduced the error due to the approximation of the transport operator by the particles. In this part, we would like to extend this concept to anynumerical approximation of transport operator, such as: SPD, Central finite difference scheme for conservation law etc.

We discuss here the error introduced by the use of the SPD approximation in our framework. The numerical analysis of this method is now well known see e.g.  \cite{jourdainreygner16}. 
\\
\begin{Lemma}\label{le transport op case one} Under Assumption \ref{ass restrict regul}  and in \textsc{CASE 1}, the following holds
	\begin{align}
		\cE^{T}_n(p) \leq C \frac{\mathfrak{h}}{M_n}\;, \, p \in \cP_n\,.
	\end{align}
\end{Lemma}

%\begin{bluetext}
\proof
By Definition \ref{de art scheme case 1}, for $p \in \cP_n$ and $n \le N-2$, we have
\begin{align}\label{eq spd start case 1}
	\cE^{T}_n(p) =  &\int \Big|\mathfrak{T}^{M_n}_{\mathfrak{h}}(\bar{u}^{N,M}_{n}(p,\cdot),p) - \tilde{\cT}_{\mathfrak{h}}(\bar{u}^{N,M}_{n}(p,\cdot),p)\Big|\ud e
% \\ &  \leq L\frac{h}{M}  + ||\mathbb{E}[u^{h,M}(t_{n+1},\tilde{P}_{t_{n+1}}^{t_{n},p},\cdot)] -  \mathbb{E}[u^{h,M}(t_{n+1},\tilde{P}_{t_{n+1}}^{t_{n},p},\cdot)]_{approx}||_{L^1} 
\end{align}
with $\bar{u}^{N,M}_n(p,\cdot) \in \cI_{M_n}$. For fixed $p$, the property of $\mathfrak{T}^{M}_{\mathfrak{h}}$ are discussed in \cite[Section 1.1]{jourdainreygner16}. We use here the estimate given in \cite[Theorem 3.1]{jourdainreygner16}, observing that, in our context, there is no error due to the discretization of the initial condition. Indeed, $\bar{u}^{N,M}_{n}(p,\cdot)$ appears already as the CDF of an empirical distribution. Thus, straightforwardly, we get
\begin{align}
\cE^{T}_n(p) \le C\frac{\mathfrak{h}}{M_n},
\end{align}
where $C$ does not depend on $p$ but depends on the Lipschitz constant of $\partial_y \mu$, recall Assumption \ref{ass restrict regul}.
\\
For the step $n=N-1$, we have,
\begin{align}
\cE^{\mathrm{T}}_{N-1}(p) = \int |\mathfrak{T}^{M}_{\mathfrak{h}}(p,\phi(\cdot))(e) - \cT_{\mathfrak{h}}(p,\phi(\cdot)(e))| \ud e
\end{align}
recall \eqref{eq last step transport approx}. It turns out that for our application $\phi$ is trivially the CDF of an empirical distribution,
recall \eqref{eq very special spec phi}. Thus, we indeed have, $\cE^{T}_{N-1}(p) \le C\frac{\mathfrak{h}}{M_{N-1}}$. 
\eproof

\subsection{Regularization error}

%As a direct consequence of the previous lemma, we obtain the following result.
\begin{Proposition}\label{pr error regu}
Under Assumption \ref{ass restrict regul}, the following holds, for $n \le N$,
\begin{align*}
 {\cE}^r_n(p) \le C  \epsilon\,,\; p \in \R^d.
\end{align*}
%\textcolor{red}{where C depends on variable $p$.}
\end{Proposition}

%\begin{bluetext}
\proof We first observe that ${\cE}^r_N(p) = 0$ as $\phi$ does not depend on $p$. 
For $n < N$, we have that
\begin{align}
 {\cE}^r_n(p) &= \int |\cV^{\epsilon}(t_n,p,e)-\cV(t_n,p,e)|\ud e  
 \\
 &= \int |\int \set{\cV(t_n,p+q,e)-\cV(t_n,p,e)} \varphi_\epsilon(q) \ud q|\ud e
 \\
 &\le \int \int | \cV(t_n,p+q,e)-\cV(t_n,p,e) |\ud e \varphi_\epsilon(q) \ud q.
\end{align}
Thus, using \eqref{eq lipschitz L1}, we compute
\begin{align}
 {\cE}^r_0(p)  \le C \int |q| \varphi_\epsilon(q) \ud q \le C\epsilon \,,
\end{align} 
%
%\begin{align}
%	{\cE}^r_0(p) & \le \int \int_{0}^{1} C(1+|p+\lambda q|)|q|\varphi_{\epsilon}(q)\ud \lambda \ud q
%\end{align}
%Then we compute
%\begin{align*}
% {\cE}^r_0(p) & \le 
% C(1+|p|) \int |q| \varphi_\epsilon(q) + \int \int_{0}^{1} C\lambda |q|^2\varphi_{\epsilon}(q)\ud \lambda \ud q \\& \le C(1+|p|) \int |q| \varphi_\epsilon(q) + C\int |q|^2\varphi_{\epsilon}(q) \ud q 
% \\ &\le C(1+|p|) \epsilon
%\end{align*}
which concludes the proof. 
\eproof
%\end{bluetext}

\subsection{Splitting error}

Recall that the splitting error is given by 
\begin{align}
\cE^{S}_n(p)= \int \Big|\tilde{\cT}_h(p,\mathbb{E}[\cV^{\epsilon}(t_{n+1},{P}_{t_{n+1}}^{t_{n},p},\cdot)])(e)- \cV^{\epsilon}(t_{n},p,e)\Big|\ud e\,. \label{eq recall error def split}
\end{align}
for $0 \le n \le N-1$, $p \in \R^d$. 
\\
In \cite{chassagneux2022numerical}, the error due the theoretical splitting has already been studied and the results obtained there can be used in our setting. However, we should  point out that here $\cV^\epsilon$ appears instead of $\cV$.
We thus first observe the following.
\begin{Lemma}\label{le basic split of splitting error}
Under Assumption \ref{ass restrict regul}, the following holds, for $p \in \R^d$,
\begin{align}\label{eq res interm}
\cE^{S}_n(p) \le C(1+ |p|^2) h^{\frac32} + \mathfrak{E}_n(p)
\end{align}
with
\begin{align} \label{eq de error n}
\mathfrak{E}_n(p) := \int \Big|\Theta_h(\cV^{\epsilon}(t_{n+1},\cdot))(p,e)- \cV^{\epsilon}(t_{n},p,e)\Big|\ud e\,.
\end{align}
\end{Lemma}

%\begin{bluetext}
\proof
We compute
\begin{align*}
\cE^{S}_n(p) \le& 
\int \Big|\tilde{\cT}_h(p,\mathbb{E}[\cV^{\epsilon}(t_{n+1},{P}_{t_{n+1}}^{t_{n},p},\cdot)])(e)- \Theta_h(\cV^{\epsilon}(t_{n+1},\cdot))(p,e)\Big|\ud e\,
\\
&+
\int \Big|\Theta_h(\cV^{\epsilon}(t_{n+1},\cdot))(p,e)- \cV^{\epsilon}(t_{n},p,e)\Big|\ud e\,.  
\end{align*}
The first term in the RHS above  is then controlled by invoking Proposition \ref{pr truncation error}.
%, namely for $0 \le n \le N-1$,
%\begin{align}\label{eq control first A1n}
%\int \Big|\tilde{\cT}_h(p,\mathbb{E}[\cV^{\epsilon}(t_{n+1},{P}_{t_{n+1}}^{t_{n},p},\cdot)])- \Theta_h(\cV^{\epsilon}(t_{n+1},\cdot))(p,e)\Big|\ud e \le C(1+ |p|^2) h^{\frac32}\;.
%\end{align}
\eproof
%\end{bluetext}

\vspace{2mm}
\noindent It remains to study $\mathfrak{E}_n$, recall \eqref{eq de error n}. The upper bound for this term is obtained in Proposition \ref{pr main point splitting} below, which allows to conclude for the splitting error in Corollary \ref{co error splitting}. We first need the following result. 

\begin{Lemma}\label{le control sup V}
Let Assumption \ref{ass restrict regul} hold. Then, 
for $n \le N-2$,
\begin{align}\label{eq majo en sup}
\sup_{(t,p,e) \in [t_n,t_{n+1}]\times \R^d \times \R }|\Theta_{t_{n+1}-t}(\cV^{\epsilon}(t_{n+1},\cdot))(p,e)- \cV^{\epsilon}(t,p,e)| \le C\epsilon \;,
\end{align}
where $C$ is a positive constant.
%\textcolor{red}{ECRIRE pour toute date t !!!}
\end{Lemma}
%\begin{bluetext}
\proof %\textcolor{red}{the following proof works only for $n<N-1$ because then the terminal condition is satisfied in the usual sense...}
Without loss of generality, we do the proof for $n=0$, working on $[0,t_1 = \mathfrak{h}]$. We denote 
$w(s,\cdot) =  \Theta_{t_1-s}(\cV^{\epsilon}(t_1,\cdot))$ for $s\in[0,t_1]$. 
Recall that the following FBSDE is well posed: for $s\in[0,t_1]$,
\begin{equation}
\left  \{
\begin{split}
P_s &= \mathrm{p} + \int_0^s b(P_t) \ud t + \bar{\sigma} W_s ,
\\
E_s &= \mathrm{e} + \int_0^s\mu(Y_t,P_t) \ud t  ,
\\
Y_s &= \cV^{\epsilon}(t_{1},P_{t_1},E_{t_1} ) - \int_s^{t_1} Z_t \ud W_t,
\end{split}
\right .
\end{equation}
with $Y_s = w(s,P_s,E_s)$ for $s \in [0,t_1]$. Without loss of generality, we will prove \eqref{eq majo en sup} at the point $(0,\mathrm{p},\mathrm{e})$ (and for $n=0$).
%(\textcolor{red}{attention a ce qui se passe en $h$)}
\\
Let $V_t = Y_t - \cV^{\epsilon}(t,P_t,E_t)$. Applying Ito's formula and using the martingale property of $Y$ we get,
\begin{align}
&\ud V_t = \ud  \cM_t  - \mu(Y_t,P_t)\partial_e \cV^{\epsilon}(t,P_t,E_t) \ud t
\\
&- 
\set{\partial_t \cV^{\epsilon}(t,P_t,E_t) +b(P_t)\partial_p \cV^{\epsilon}(t,P_t,E_t)+ \frac12 \Delta_p \cV^{\epsilon}(t,P_t,E_t) } \ud t
\end{align}
where    $\cM$ is a square integrable martingale with $\cM_0 = 0$. 
\\
Using the PDE \eqref{eq pde regul regul V} satisfied by $\cV^{\epsilon}$ we get
\begin{align}\label{eq for V}
&\ud V_t = \ud  \cM_t  -
%\\
%&+
 \!\! \left( \set{ 
\mu(Y_t,P_t) ) - \mu(\cV^{\epsilon}(t,P_t,E_t),P_t) )
}
\partial_e \cV^{\epsilon}(t,P_t,E_t) 
+ \theta^\epsilon(t,P_t,E_t)  \right ) \ud t.
\end{align}
Observe that
\begin{align}
\mu(Y_t,P_t) - \mu(\cV^{\epsilon}(t,P_t,E_t),P_t) )  = c_t V_t  \; \text{ with } \;c_t :=  \int_0^1 \partial_y \mu(\cV^{\epsilon}(t,P_t,E_t) + \lambda V_t,P_t) \ud \lambda
\end{align}
%with
%\begin{align}
%c_t :=  \int_0^1 \partial_y \mu(Y_t,P_t) \ud \lambda
%\\
%\frac{\mu(Y_t,P_t)  - \mu(\cV^{\epsilon}(t,P_t,E_t),P_t) )}{Y_t-\cV^{\epsilon}(t,P_t,E_t)} \1_{\set{Y_t \neq \cV^{\epsilon}(t,P_t,E_t)}}
%\end{align}
and from the property of $\mu$, recall \eqref{eq conservation law}, we have
\begin{align}\label{eq property ct}
{c}_t  \le -\ell_1 < 0\;.
\end{align}
We get, for $s \in [0,t_1]$
\begin{align}\label{eq key rep for V}
V_s - V_0
=
- \int_0^s \left( {c}_t  V_t \partial_e \cV^{\epsilon}(t,P_t,E_t) 
+ \theta^\epsilon(t,P_t,E_t)\right) \ud t  + \cM_s.
\end{align}
We set, for $0 \le t \le t_1$, $\cE_t = e^{\int_0^t {c}_s \partial_e \cV^{\epsilon}(s,P_s,E_s)   \ud s}$ and, we have
\begin{align}
0 \le \cE_t \le 1\;,\; \text{ for all } 0 \le t \le t_1\,,
\end{align}
since $\partial_e \cV^{\epsilon} $ is non negative, recall also \eqref{eq property ct}.
We then compute
\begin{align}
\cE_s V_s - V_0
=
-\int_0^s  
\theta^\epsilon(t,P_t,E_t) \cE_t \ud t  + \cN_s
\end{align}
with $\cN$ a square integrable martingale such that $\cN_0 = 0$. In particular, recalling \eqref{eq control theta}, we get
\begin{align}
V_0 \le |\esp{\cE_s V_s}| + C\mathfrak{h}\epsilon +  C \epsilon \esp{\int_0^s  
\partial_e \cV^\epsilon(t,P_t,E_t) \cE_t \ud t}
\end{align}
where we use the fact that $|\partial_p \cV|$ is bounded and $\partial_e \cV^\epsilon(t,P_t,E_t) \ge 0$.
The above inequality reads also
\begin{align}
|V_0| \le |\esp{\cE_s V_s}| + C\mathfrak{h}\epsilon - C \epsilon \esp{\int_0^s  \frac{{c}_t}{|{c}_t|} \partial_e \cV^\epsilon(t,P_t,E_t) \cE_t \ud t}.
\end{align}
Observing that 
$$ \dot{\cE}_t = {c}_t
\partial_e \cV^\epsilon(t,P_t,E_t) \cE_t ,$$
since $|c_t| \ge \ell_1 >0$, we obtain
\begin{align}
|V_0| \le |\esp{\cE_s V_s}| + C(1+\mathfrak{h})\epsilon  
\end{align}
and {since $V_{t_1} = 0$}
$$|w(0,\mathrm{p},\mathrm{e})-\cV^{\epsilon}(0,\mathrm{p},\mathrm{e})| = |V_0| \le C\epsilon\,,$$
which concludes the proof.
\eproof
%\end{bluetext}

\begin{Proposition}\label{pr main point splitting} Let Assumption \ref{ass restrict regul} hold. Then, 
for $n \le N-2$, $p\in \R^d$,
\begin{align}
\mathfrak{E}_n(p)=\int |\Theta_{\mathfrak{h}}(\cV^{\epsilon}(t_{n+1},\cdot))(p,e)- \cV^{\epsilon}(t_{n},p,e)| \ud e \le  C \mathfrak{h} \epsilon \,.
\end{align}
\end{Proposition}

%\begin{bluetext}
\proof  Without loss of generality, we  prove the statement in the same setting as the one used for the proof of Lemma \ref{le control sup V}. We just stress here the dependence upon the initial condition:
$V^{\mathrm{e}}_t = Y^{\mathrm{e}}_t - \cV^{\epsilon}(t,P_t,E^{\mathrm{e}}_t)$ since $E_0 = {\mathrm{e}}$.
We consider the tangent process $\partial_e {E}$ given by
\begin{align}
\partial_e {E}_t &= 1 + \int_0^t \partial_y \mu(Y_s,P_s)\partial_e w(s,P_s,E_s) \partial_e {E}_s \ud s ,
\\
&= e^{\int_0^t\partial_y \mu(Y_s,P_s)\partial_e w(s,P_s,E_s)  \ud s}.
\end{align}
And we observe that $0 \le \partial_e E_t \le 1$, for all $0 \le t \le t_1$, since $\partial_e w \ge 0$ and  $\mu$ is decreasing in $y$, recall \eqref{eq conservation law}. \\
In order to bound the error $\int |V^e_0| \ud e$, we will study the dynamics of
$t \mapsto \int |\esp{V^{e}_t\partial_e {E}^{e}_t} | \ud e$.
%%%

%\begin{align}%{eq key rep for V}
%V_s - V_0
%=
%\int_0^s \left( \tilde{c}_t  V_t \partial_e \cV^{\epsilon}(t,P_t,E_t) 
%+ \theta^\epsilon(t,P_t,E_t)\right) \ud t  + M_s
%\end{align}

Recalling 
\eqref{eq key rep for V}, we compute
\begin{align}
&\ud ( V^e_t \partial_e E_t) = \ud  \cN_t  + V^e_t \partial_y \mu(Y_t,P_t)\partial_e w(t,P_t,E_t) \partial_e {E}_t \ud t
\\
&-
 \!\! \left( {c}_t  V^e_t \partial_e \cV^{\epsilon}(t,P_t,E_t)  \partial_e {E}_t
+ \theta^\epsilon(t,P_t,E^e_t) \partial_e {E}_t \right ) \ud t,
\end{align}
where $\cN$ is a square integrable martingale satisfying $\cN_0 = 0$.
\begin{align}
\esp{V^e_{\mathfrak{h}} \partial_e E_{\mathfrak{h}}- V^e_0} &= \esp{\int_0^{\mathfrak{h}} V^e_t \big(  \partial_y \mu(Y_t,P_t)\partial_e[ w(t,P_t,E^e_t) ] -{c}_t  V^e_t \partial_e [\cV^{\epsilon}(t,P_t,E^e_t)] \big)
\ud t}
\\
&+\esp{ \int_0^{\mathfrak{h}} \theta^\epsilon(t,P_t,E_t) \partial_e {E}_t  \ud t}\,.
\end{align}
Since $V_{\mathfrak{h}} = 0$ and $\partial_e w \ge 0$,  $\partial_e \cV^\epsilon \ge 0$, $\partial_e {E}_t \ge 0$, we obtain
\begin{align}
\int |V^e_0| \ud e &\le 
\int_0^{\mathfrak{h}} \esp{\int |V^e_t  \partial_y \mu(Y_t,P_t)|\partial_e[ w(t,P_t,E^e_t) ]  \ud e } \ud t \label{eq main decomp here}
\\
&+\int_0^{\mathfrak{h}}  \esp{\int |{c}_t  V^e_t| \partial_e [\cV^{\epsilon}(t,P_t,E^e_t)] \ud e }\ud t \nonumber
\\
&+\int_0^{\mathfrak{h}} \esp{ \int |\theta^\epsilon(t,P_t,E_t)| \partial_e {E}_t  \ud e } \ud t . \nonumber
\end{align}
We obtain, using Lemma \ref{le control sup V} and the bound on $\|\partial_y \mu\|_\infty$,
\begin{align}
\int |V^e_t  \partial_y \mu(Y_t,P_t)| \partial_e[ w(t,P_t,E^e_t) ]  \ud e 
&\le C \epsilon \int \partial_e[ w(t,P_t,E^e_t) ]  \ud e, \nonumber
\\
&\le C \epsilon, \label{eq pre estim 1}
\end{align}
where we used the fact that $w$ is bounded in the last inequality. Similar arguments, since ${c}$ is bounded, lead to
\begin{align}
\int |{c}_t  V^e_t| \partial_e [\cV^{\epsilon}(t,P_t,E^e_t)] \ud e \le C \epsilon \,. \label{eq pre estim 2}
\end{align}
From \eqref{eq control theta}, we know that
\begin{align}%\label{eq control theta}
|\theta^\epsilon(t,P_t,E_t)|
\le C\epsilon\left (\int |\partial_p \cV(t,P_t+q,E^e_t)|\varphi_\epsilon(q) \ud q  + \partial_e\cV^{\epsilon}(t,P_t,E_t)\right ).
\end{align}
We compute, using the change of variable $e \mapsto \zeta = {E}^e_t$
\begin{align}
\int \int |\partial_p \cV(t,P_t+q,E^e_t)|\varphi_\epsilon(q) \ud q \partial_e {E}^e_t  \ud e 
&= \int  \int |\partial_p \cV(t,P_t+q,\zeta)|\varphi_\epsilon(q) \ud q  \ud \zeta,
\\
&\le C,
\end{align}
where we used the $L^1$-bound on $\partial_p V$ given in \eqref{eq control deriv L1} .
We then compute
\begin{align}
\esp{ \int |\theta^\epsilon(t,P_t,E_t)| \partial_e {E}_t  \ud e }
&\le C\epsilon(1+ \esp{\int  \partial_e\cV^{\epsilon}(t,P_t,E_t)\partial_e {E}_t  \ud e} ),
\\
&\le C\epsilon.
\end{align}
The proof is concluded by combining  the previous inequality, the estimates given in \eqref{eq pre estim 1}-\eqref{eq pre estim 2} with \eqref{eq main decomp here}.
\eproof
%\end{bluetext}

\vspace{2mm}
\noindent Combining Lemma \ref{le basic split of splitting error} with the result Proposition \ref{pr main point splitting}, we obtain straightforwardly the following.
\begin{Corollary}\label{co error splitting}
 Let Assumption \ref{ass restrict regul} hold. Then, 
for $n \le N-2$, $p \in \R^d$,
\begin{align}\label{eq control error splitting}
\cE^{S}_n(p) \le C(1+ |p|^2) \mathfrak{h}^{\frac32} + C\epsilon \mathfrak{h} \,.
\end{align}
\end{Corollary}

\subsection{Diffusion error}
We now study the term given in \eqref{eq error def diff} and we straightforwardly observe
\begin{align}
\cE^{D}_n(p) &= \int \Big|\tilde{\cT}_h(p,\mathbb{E}[\cV^{\epsilon}(t_{n+1},\widehat{P}_{t_{n+1}}^{t_{n},p},\cdot)])(e) - \tilde{\cT}_h(p,\mathbb{E}[\cV^{\epsilon}(t_{n+1},{P}_{t_{n+1}}^{t_{n},p},\cdot)])(e) \Big|\ud e\,, 
%
%\cE^{D}_n(p) &= \int \Big|\cT_h(\mathbb{E}[u^{\epsilon}(t_{n+1},\tilde{P}_{t_{n+1}}^{t_{n},p},\cdot)],p) - \cT_h(\mathbb{E}[u^{\epsilon}(t_{n+1},{P}_{t_{n+1}}^{t_{n},p},\cdot)],p) \Big|\ud e
\\
&\le
\int \Big| \mathbb{E}[\cV^{\epsilon}(t_{n+1},\widehat{P}_{t_{n+1}}^{t_{n},p},e)] - \mathbb{E}[\cV^{\epsilon}(t_{n+1},{P}_{t_{n+1}}^{t_{n},p},e)] \Big|\ud e\,,
\label{eq expression diff error bis}
\end{align}
where we use the $L^1$-stability of $\tilde{\cT}_h$.% \textcolor{magenta}{(uniform in p?)}.

\vspace{2mm}
\noindent This motivates the introduction of the following auxiliary result.

\begin{Proposition}\label{pr reg auxiliary}
Let Assumption \ref{ass restrict regul} hold. 
For $0 \le n \le N-1$, let $w_n$ be the solution on $[t_n,t_{n+1}]\times \R^d$ to the following PDE:
	\begin{align}\label{w pde}
	\partial_t w + b(p)\partial_p w + \frac12\bar{\sigma}^2 \Delta_{pp} w = 0 \text{ and } w(t_{n+1},\cdot):=\cV^{\epsilon}(t_{n+1},\cdot,e)\;,
	\end{align}
	for a fixed $e \in \R$.
	Then, under Assumption \ref{ass restrict regul}, the followings hold, for $q,q' \in \R^d$ and $t\in [0,T)$,
	\begin{align}
		\int |\partial_{p^\ell}w_n(t,q,e)-\partial_{p^\ell}w_n(t,q',e)| \ud e \leq \frac{C}{\epsilon}|q-q'|
	\end{align}
	and
	\begin{align}
		\int |\partial^2_{p^\ell p^{\ell'}}w_n(t,q,e)-\partial^2_{p^\ell p^{\ell'}}w_n(t,q',e)| \ud e  \leq \frac{C}{\epsilon^2}|q-q'|,
	\end{align}
for all $\ell,\ell' \in \set{1,\dots,d} $.
Moreover,
	\begin{align}\label{eq bound partial p w}
		\int |\partial_p w_n(t,p,e)|\ud e \leq C.
		%\;\text{ and }\; \mathbb{E}\int |\partial_{pp}^2 w(t,\tilde{P}_{t}^{t_{n},p},e)|\ud e \leq \frac{C}{\epsilon}.
	\end{align}
\end{Proposition}

%\begin{bluetext}
\proof
For ease of presentation, the proof is done in the one-dimensional case.\\
By Theorem 2.3.5 in Zhang \cite{zhang2001some}, we have the following expressions for the first and second order derivative with respect to vector $p$, for $(t,p,e) \in [t_n,t_{n+1}]\times \R \times \R$,
\begin{align}
	\partial_p w_n(t,p,e) &= \esp{\partial_p \cV^{\epsilon}(t_{n+1},P^{t,p}_{t_{n+1}},e)\partial_p P^{t,p}_{t_{n+1}}},
	\label{eq rep p deriv w}
	\\
	\partial^2_{pp} w_n(t,p,e) &= \esp{\partial^2_{pp} \cV^{\epsilon}(t_{n+1},P^{t,p}_{t_{n+1}},e)(\partial_p P^{t,p}_{t_{n+1}})^2 + \partial_p \cV^{\epsilon}(t_{n+1},P^{t,p}_{t_{n+1}},e)\partial_{pp}^2 P^{t,p}_{t_{n+1}} },
	\label{eq rep second p deriv w}
	\end{align}
where, for $t\le s \le t_{n+1}$,
\begin{align}
P^{t,p}_s &= p + \int_{t}^s b(P^{t,p}_r) \ud r + \bar{\sigma}(W_s-W_t)\,,
\\
\partial_p P^{t,p}_s &= 1 + \int_{t}^s b'(P^{t,p}_r) \partial_p P^{t,p}_r\ud r \,,
\\
\partial^2_{pp} P^{t,p}_s &= \int_{t}^s \left( b''(P_r)|\partial_p P^{t,p}_r|^2 + b'(P_r)\partial^2_{pp} P^{t,p}_r \right)\ud r\;.
\end{align}
In this setting, classical arguments lead to,
\begin{align}\label{eq control second order}
\esp{\sup_{s\in [t,t_{n+1}]} |\partial_{p} P^{t,p}_s|^\kappa + \sup_{s\in [t,t_{n+1}]} |\partial^2_{pp} P^{t,p}_s|^\kappa} \le C_\kappa,
\end{align}
%recalling also \eqref{eq bound partial P} for future use.\\
 and for $q,q' \in \R$, $s \in [t,T]$, $\kappa \ge 1$,
\begin{align}\label{eq general control}
\esp{|P^{t,q}_s - P^{t,q'}_s |^\kappa + |\partial_pP^{t,q}_s - \partial_pP^{t,q'}_s |^\kappa  + |\partial^2_{pp}P^{t,q}_s - \partial^2_{pp}P^{t,q'}_s |^\kappa } \le C|q-q'|^\kappa\;.
\end{align}

\noindent We now study the second order derivative given in \eqref{eq rep second p deriv w}. 
\\We compute, for $(t,q,q')\in[t_n,t_{n+1}]\times \R \times \R$,
\begin{align}
&\int |\partial^2_{pp}w_n(t,q,e)-\partial^2_{pp}w_n(t,q',e)| \ud e
\\ &\le
\esp{\int|\partial^2_{pp} \cV^\epsilon(t_{n+1},P^{t,q}_{t_{n+1}},e)(\partial_p P^{t,q}_{t_{n+1}})^2
-
\partial^2_{pp} \cV^\epsilon(t_{n+1},P^{t,q'}_{t_{n+1}},e)(\partial_p P^{t,q'}_{t_{n+1}})^2| \ud e}:= B_1
\\
&+ 
\esp{\int|\partial_{p} \cV^\epsilon(t_{n+1},P^{t,q}_{t_{n+1}},e)\partial_{pp}^2 P^{t,q}_{t_{n+1}}
-
\partial_{p} \cV^\epsilon(t_{n+1},P^{t,q'}_{t_{n+1}},e)\partial_{pp}^2 P^{t,q'}_{t_{n+1}}| \ud e}:= B_2 \label{eq de B2}
\end{align}

For the first term $B_1$,
\begin{align*}
	B_1 &\leq \esp{\int|\partial^2_{pp} \cV^\epsilon(t_{n+1},P^{t,q}_{t_{n+1}},e)(\partial_p P^{t,q}_{t_{n+1}})^2
	-
	\partial^2_{pp} \cV^\epsilon(t_{n+1},P^{t,q}_{t_{n+1}},e)(\partial_p P^{t,q'}_{t_{n+1}})^2| \ud e}
	\\&  
	+ \esp{\int|\partial^2_{pp} \cV^\epsilon(t_{n+1},P^{t,q}_{t_{n+1}},e)(\partial_p P^{t,q'}_{t_{n+1}})^2
	-
	\partial^2_{pp} \cV^\epsilon(t_{n+1},P^{t,q'}_{t_{n+1}},e)(\partial_p P^{t,q'}_{t_{n+1}})^2| \ud e}.
\end{align*}
We then compute
\begin{align}
	&\esp{\int|\partial^2_{pp} \cV^\epsilon(t_{n+1},P^{t,q}_{t_{n+1}},e)(\partial_p P^{t,q}_{t_{n+1}})^2
	-
	\partial^2_{pp} \cV^\epsilon(t_{n+1},P^{t,q}_{t_{n+1}},e)(\partial_p P^{t,q'}_{t_{n+1}})^2| \ud e},
	 \\& \leq  
	\esp{\Big|(\partial_p P^{t,q}_{t_{n+1}})^2 - (\partial_p P^{t,q'}_{t_{n+1}})^2\Big|\int|\partial^2_{pp} \cV^\epsilon(t_{n+1},P^{t,q}_{t_{n+1}},e)| \ud e} ,
\\ & \leq \frac{C}{\epsilon}\esp{\Big|(\partial_p P^{t,q}_{t_{n+1}})^2 - (\partial_p P^{t,q'}_{t_{n+1}})^2\Big|}, \label{eq key last row}
\end{align}
where we used for the last inequality Lemma \ref{le bound for V epsilon}.
Observing that
\begin{align*}
	\Big|(\partial_p P^{t,q}_{t_{n+1}})^2 - (\partial_p P^{t,q'}_{t_{n+1}})^2\Big| &\leq \Big|\partial_p P^{t,q}_{t_{n+1}} +\partial_p P^{t,q'}_{t_{n+1}}\Big|\cdot\Big|\partial_p P^{t,q}_{t_{n+1}} - \partial_p P^{t,q'}_{t_{n+1}}\Big|,
	%\\& \le \sup_{s\in[t,t_{n+1}]}(|\partial_p P^{t,q'}_{s}|+|\partial_p P^{t,q}_{s}|)\Big|\partial_p P^{t,q}_{t_{n+1}} - \partial_p P^{t,q'}_{t_{n+1}}\Big|
\end{align*}
we combine Cauchy-Schwarz inequality with \eqref{eq general control} to obtain
\begin{align}
\esp{\Big|(\partial_p P^{t,q}_{t_{n+1}})^2 - (\partial_p P^{t,q'}_{t_{n+1}})^2\Big| } \le C|q-q'|\,.
\end{align}
Combining the previous inequality with  \eqref{eq key last row}, we get
\begin{align}\label{eq first step B1}
\esp{\int|\partial^2_{pp} \cV^\epsilon(t_{n+1},P^{t,q}_{t_{n+1}},e)(\partial_p P^{t,q}_{t_{n+1}})^2
	-
	\partial^2_{pp} \cV^\epsilon(t_{n+1},P^{t,q}_{t_{n+1}},e)(\partial_p P^{t,q'}_{t_{n+1}})^2| \ud e} \le \frac{C}{\epsilon}|q-q'|\,.
\end{align}
%\textcolor{red}{bound on$\int|\partial^2_{pp} \cV^\epsilon(t_{n+1},P^{t,q}_{t_{n+1}},e)| \ud e$ }

%Since we have \eqref{estimate order 2}, and besides we also have $\mathbb{E}\Big|(\partial_p P^{t,q}_{t_{n+1}})^2 - (\partial_p P^{t,q'}_{t_{n+1}})^2\Big|\le C|b-a|$
%\begin{align*}
%	\Big|(\partial_p P^{t,q}_{t_{n+1}})^2 - (\partial_p P^{t,q'}_{t_{n+1}})^2\Big| &\leq \Big|\partial_p P^{t,q}_{t_{n+1}} +\partial_p P^{t,q'}_{t_{n+1}}\Big|\cdot\Big|\partial_p P^{t,q}_{t_{n+1}} - \partial_p P^{t,q'}_{t_{n+1}}\Big|\\& \le \sup_{s\in[t,t_{n+1}]}(|\partial_p P^{t,q'}_{s}|+|\partial_p P^{t,q}_{s}|)\Big|\partial_p P^{t,q}_{t_{n+1}} - \partial_p P^{t,q'}_{t_{n+1}}\Big|
%\end{align*}
We also compute
\begin{align*}
	&\esp{\int|\partial^2_{pp} \cV^\epsilon(t_{n+1},P^{t,q}_{t_{n+1}},e)(\partial_p P^{t,q'}_{t_{n+1}})^2
	-
	\partial^2_{pp} \cV^\epsilon(t_{n+1},P^{t,q'}_{t_{n+1}},e)(\partial_p P^{t,q'}_{t_{n+1}})^2| \ud e} \\& \leq  
	\esp{(\partial_p P^{t,q'}_{t_{n+1}})^2\int|\partial^2_{pp} \cV^\epsilon(t_{n+1},P^{t,q}_{t_{n+1}},e)
	-
	\partial^2_{pp} \cV^\epsilon(t_{n
	+1},P^{t,q'}_{t_{n+1}},e)| \ud e} 
	\\
	&\le \frac{C}{\epsilon^2}|q-q'|
\end{align*}
where we used Lemma \ref{le bound for V epsilon}	and \eqref{eq control second order}, for the last inequality.
Combining \eqref{eq first step B1} with the previous inequality, we finally obtain that
\begin{align}\label{eq conclu B1}
	B_1 \leq \frac{C}{\epsilon^2}|q-q'| +\frac{C}{\epsilon}|q-q'|\;.
\end{align}
For the term $B_2$ given in \eqref{eq de B2}, we observe
\begin{align*}
B_2 & \leq \esp{\int|\partial_{p} \cV^\epsilon(t_{n+1},P^{t,q}_{t_{n+1}},e)\partial_{pp}^2 P^{t,q}_{t_{n+1}}
-
\partial_{p} \cV^\epsilon(t_{n+1},P^{t,q}_{t_{n+1}},e)\partial_{pp}^2 P^{t,q'}_{t_{n+1}}| \ud e}
\\& + \esp{\int|\partial_{p} \cV^\epsilon(t_{n+1},P^{t,q}_{t_{n+1}},e)\partial_{pp}^2 P^{t,q'}_{t_{n+1}}
-\partial_{p} \cV^\epsilon(t_{n+1},P^{t,q'}_{t_{n+1}},e)\partial_{pp}^2 P^{t,q'}_{t_{n+1}}| \ud e}
\end{align*}
Combining Lemma \ref{le bound for V epsilon}, the estimates given in \eqref{eq general control}-\eqref{eq control second order},
with similar computations as done previously, we obtain
\begin{align}\label{eq conclu B2}
	B_2 \leq \frac{C}{\epsilon}|q-q'| + C|q-q'|\;.
\end{align}
Inserting  \eqref{eq conclu B1}-\eqref{eq conclu B2} back into \eqref{eq de B2}, we obtain the
\begin{align}
\int |\partial^2_{pp}w_n(t,q,e)-\partial^2_{pp}w_n(t,q',e)| \ud e \le \frac{C}\epsilon |q-q'|\,,
\end{align}
as $\epsilon \in (0, 1)$.
\\
Similar arguments allow to retrieve the estimates for the first order derivatives.
\eproof
%\end{bluetext}

\begin{Proposition}\label{pr d dim}
Under Assumption \ref{ass restrict regul},
the following holds
\begin{align*}
\cE^{D}_n(p) &\le C\frac{\mathfrak{h}^{\frac32}}{\epsilon^2},
\end{align*}
for $0 \le n \le N-1$ and $p \in \R^d$.
\end{Proposition}

%\begin{bluetext}
\proof  For ease of presentation, we do the proof in dimension $d=1$.\\
1. We first recall  Definition \ref{de bro inc} and \eqref{eq de one step disc euler scheme P}. In this context, we introduce two new discrete processes:
\begin{align*}
	&\widehat{P}^{t_n,p}_{t}= p +  b(p)(t-t_{n}) + \bar{\sigma}\frac{\sqrt{t-t_{n}}}{\sqrt{\mathfrak{h}}}\Delta \widehat{W}_n\,, \\
	&\bar{P}_{t}^{t_n,p,\lambda} = p  + b(p )(t-t_{n}) + \lambda \bar{\sigma}\frac{\sqrt{t-t_{n}}}{\sqrt{\mathfrak{h}}}\Delta \widehat{W}_n,
\end{align*}
for $p \in \R$, $t \in [t_n,t_{n+1}]$, $\lambda \in [0,1]$ and $0 \le n \le N-1$.\\
Now, recalling \eqref{eq expression diff error bis}, we first observe
\begin{align}\label{eq classic trick}
\mathbb{E}[\cV^{\epsilon}(t_{n+1},\widehat{P}_{t_{n+1}}^{t_{n},p},e)] - \mathbb{E}[\cV^{\epsilon}(t_{n+1},{P}_{t_{n+1}}^{t_{n},p},e)]
=
\mathbb{E}[w_n(t_{n+1},\widehat{P}_{t_{n+1}}^{t_{n},p},e)] - w_n(t_{n},p,e)\;.
\end{align}
We apply the discrete Ito formula  given of Proposition 14  in \cite{chassagneux2019numerical} to our much simpler framework, to obtain
\begin{align*}
\mathbb{E}[w_n(t_{n+1},\widehat{P}_{t_{n+1}}^{t_{n},p},e)]  - w_n(t_{n},p,e) &= \int_{t_n}^{t_{n+1}} \esp{\partial_t w_n(t,\widehat{P}_{t}^{t_{n},p},e)
+b(p)\partial_p w_n(t,\widehat{P}_{t}^{t_{n},p},e) 
} \ud t\\
&+\frac{\bar{\sigma}^2}2\esp{\int_{t_n}^{t_{n+1}}\!\!\int_0^1\partial^2_{pp}w_n(t,\bar{P}_{t}^{t_n,p,\lambda},e)\ud \lambda \ud t }\;.
\end{align*}
Using the PDE \eqref{w pde} satisfied by $w_n$, the previous equality leads to
\begin{align}
\mathbb{E}[w_n(t_{n+1}&,\widehat{P}_{t_{n+1}}^{t_{n},p},e)]  - w_n(t_{n},p,e) \nonumber
\\ 
&=\int_{t_n}^{t_{n+1}} \esp{\left(b(p)-b(\widehat{P}_{t}^{t_{n},p})\right)\partial_p w_n(t,\widehat{P}_{t}^{t_{n},p},e) 
} \ud t \label{eq first term disc diff error}\\
&+\frac{\bar{\sigma}^2}2\esp{\int_{t_n}^{t_{n+1}}\!\!\int_0^1[\partial^2_{pp}w_n(t,\bar{P}_{t}^{t_n,p,\lambda},e)-\partial^2_{pp}w_n(t,\widehat{P}_{t}^{t_n,p},e)]\ud \lambda \ud t } \label{eq second term disc diff error} \;.
\end{align}
2. For  $e \in \R$, we denote by $A^1(e)$, the term in \eqref{eq first term disc diff error} and by $A^2(e)$, the term in \eqref{eq second term disc diff error} .\\
We compute
\begin{align*}
|A^1(e)| &\le \int_{t_n}^{t_{n+1}} \esp{|(b(p)-b(\widehat{P}_{t}^{t_{n},p})| |\partial_p w_n(t,\widehat{P}_{t}^{t_{n},p},e)| } \ud t,
\\
& \le C \mathfrak{h}^\frac12  \int_{t_n}^{t_{n+1}}\esp{|\partial_p w_n(t,\widehat{P}_{t}^{t_{n},p},e)| } \ud t  \,,
\end{align*}
and thus
\begin{align}
\int |A^1(e)| \ud e & \le C \mathfrak{h}^\frac12  \int_{t_n}^{t_{n+1}}\esp{\int |\partial_p w_n(t,\widehat{P}_{t}^{t_{n},p},e)| \ud e} \ud t \,,
\nonumber
\\
& \le C \mathfrak{h}^\frac32  \,,
\label{eq control error A1}
 \end{align}
 where we used \eqref{eq bound partial p w} for the last inequality.
 \\
 Using Proposition \ref{pr reg auxiliary}, we compute
 \begin{align*}
 \int |\partial^2_{pp}w_n(t,\bar{P}_{t}^{t_n,p,\lambda},e)-\partial^2_{pp}w_n(t,\widehat{P}_{t}^{t_n,p},e)| \ud e
 &\le \frac{C}{\epsilon^2}|\bar{P}_{t}^{t_n,p,\lambda}-\widehat{P}_{t}^{t_n,p}| \,,
 \\
 &\le \frac{C}{\epsilon^2} \sqrt{\mathfrak{h}} \,.
 \end{align*}
 This yields
 \begin{align}
 \int |A_2(e)| \ud e \le  \frac{C}{\epsilon^2} \mathfrak{h}^\frac32  \,.
 \end{align}
 Combining the previous inequality with \eqref{eq control error A1} and \eqref{eq first term disc diff error}-\eqref{eq second term disc diff error},
 we get
 \begin{align*}
 \int |\mathbb{E}[w_n(t_{n+1}&,\widehat{P}_{t_{n+1}}^{t_{n},p},e)]  - w_n(t_{n},p,e)| \ud e \le  \frac{C}{\epsilon^2} \mathfrak{h}^\frac32 \,.
 \end{align*}
 The proof is concluded recalling \eqref{eq classic trick} and \eqref{eq expression diff error bis}.
\eproof
%\end{bluetext}

%\input{oldstuff}

\subsection{Proof of Theorem \ref{th main conv result}}
\label{subse proof th main conv result}
%\begin{bluetext}
The proof of our main result is now almost straightforward.
From Lemma \ref{le transport op case one}, Corollary \ref{co error splitting} and Proposition \ref{pr d dim},
we observe that
\begin{align}
\esp{ \cE^{\mathrm{T}}_n(\widehat{P}_{t_n}) } \le C\frac{\mathfrak{h}}M\;,\;
\esp{  \cE^{D}_n(\widehat{P}_{t_n})} \le  C\frac{\mathfrak{h}^{\frac32}}{\epsilon^2} \;\text{ and }\;
\esp{  \cE^{S}_n(\widehat{P}_{t_n}) }  \le C\mathfrak{h}\left(\mathfrak{h}^{\frac12} + \epsilon  \right)\,.
\end{align}
Summing over $n$, we get
\begin{align}\label{eq main th interm 1}
 \sum_{n=0}^{N-2}\esp{ \cE^{\mathrm{T}}_n(\widehat{P}_{t_n})   +  \cE^{D}_n(\widehat{P}_{t_n}) +  \cE^{S}_n(\widehat{P}_{t_n}) }  
 \le
C \left(\frac1M + \frac{\sqrt{\mathfrak{h}}}{\epsilon^2} + \epsilon \right).
\end{align}
For the last step, we obtain 
\begin{align}
\esp{{\cE}_{N-1}(\widehat{P}_{t_{N-1}})} \le C \left(\frac{h}M +\epsilon+ \mathfrak{h}\left(\mathfrak{h}^{\frac12} + \epsilon  \right) \right)
\end{align}
combining  \eqref{eq control last step} with Proposition \ref{pr error regu} and Proposition \ref{pr truncation error}.
The previous inequality combined with \eqref{eq main th interm 1} and  Proposition \ref{pr global stab}
yields \eqref{eq control main error}.
\\
Balancing optimally the errors in \eqref{eq control main error} allows to conclude  by proving \eqref{eq main err in term of h}.
\eproof
%\end{bluetext}

%\input{numeric}

\section{Numerics}
\label{se numerics}
%\textcolor{red}{
%- introduire les exemples
%\\
%- tester numériquement case 1 et case 2
%\\
%- faire un exemple en multipériode?
%}

{\color{r_two}
In this section, we realize  numerical experiments to test in practice the schemes introduced in Section \ref{se algo}.

\subsection{Setting}  
We first introduce the models  we will use. They have already been considered in \cite{chassagneux2022numerical}, which facilitates the comparison with previous numerical results.
}

\noindent Let us first define the following toy model where the process $P$ corresponds to  a  Brownian motion and is a multidimensional version of the model given in Remark \ref{re example for ass}:

\begin{Example}[Linear model]\label{ex lin model p2}
	\label{ex linear emission}
	\begin{align}
		dP_t& = \sigma\ud W_{t}, dE_t= \left(\frac1{\sqrt{d}}\sum_{\ell=1}^dP^\ell_t-Y_t\right)\ud t 
	\end{align}
	with terminal function $(p,e)\mapsto \phi(p,e) = \1_{\set{e \ge 0}}$ and 
	where $W$ is a $d$-dimensional Brownian motion and $\sigma > 0$.%an invertible matrix in $\cM_d$.
\end{Example}

\noindent We will also consider a multiplicative model:
%presented in \cite{chassagneux2022numerical}:
\begin{Example}[Multiplicative model]\label{ex mulitplicative model p2} Let $W$ be a $d$-dimensional Brownian motion. For all $\ell \in \set{1,\dots,d}$, we set
	\label{ex gbm P}%{ex bm pos emit}
	\begin{align}
		\ud P^\ell_t = \mu P^\ell_t \ud t + \sigma P^\ell_t \ud W^\ell_t, \,P^\ell_0=1,  \text{ and } \ud E_t = \tilde \mu(Y_t,P_t) \ud t
	\end{align}
	with
	%\begin{align*}
	$
	(y,p) \mapsto \tilde{\mu}(y,p)= \left(\prod_{\ell=1}^d p^\ell \right)^{\frac1{\sqrt{d}}}e^{-\theta y}
	$, for some $\theta > 0$. The terminal condition is given by $(p,e)\mapsto \phi(p,e) = \1_{\set{e \ge 0}}$\;.
	%\end{align*}
\end{Example} 

{\color{r_two}
As  mentioned in \cite{chassagneux2022numerical}, one advantage of these models is that they can be reduced to models with lower dimension (at most 2) whatever the value of $d$. For comparison, we will also use some numerical \emph{proxy}: for the first model in Example \ref{ex lin model p2}, we use the method considered in \cite{bossy1996convergence}  applied to the one dimensional reduction of the model, for the second model in Example \ref{ex mulitplicative model p2} we use the NN \& Upwind scheme. Both methods are introduced and discussed in \cite{chassagneux2022numerical}.

\vspace{2mm} 
\noindent We also observe that the above models fit into the setting of Section \ref{subsubse brownian setting}, recall Assumption \ref{ass brownian functional}. As shown below and as expected, they are thus numerically tractable. We shall use in this regard the scheme introduced in Definition \ref{de art scheme case 1}, linked to \textsc{CASE 1}. We will also consider \textsc{CASE 2}, where the number of particles is capped. To this end, we need to define precisely the operator $\mathfrak{R}^{M,m}$ appearing in Definition  \ref{de art scheme case 2}, which is the purpose of the next section.
}

%In this part, we first compare the differences between schemes CASE 1 and CASE 2 defined in Section \ref{subsubse approx generic} for model Example \ref{ex lin model p2}: a linear toy model with Brownian setting. Then to validate our numerical results, we compare the numerical scheme CASE 2 to ``proxy'' solution presented in \cite{bossy1996convergence} for different levels of volatility for model Example \ref{ex lin model p2}. Then we have also tested the performance of the scheme CASE 2 in a multiplicative model Example \ref{ex mulitplicative model p2} against NN \& Upwind method proposed in the article \cite{chassagneux2022numerical}. 
%
%Before we dive into numerical part, further clarifications on scheme CASE 2 are needed: as we mentioned earlier in Section \ref{subsubse approx generic}
%, various implementations of reducing the number of particles are possible for scheme CASE 2. In the following, we will define and clarify the operator $\mathfrak{R}^{M,m}$ as well as its associated different implementations.

\subsection{Definition of $\mathfrak{R}$}

{\color{black}
There are various possible implementations  for the operator $\mathfrak{R}$ whose goal is to tame the computational cost, recall Definition  \ref{de art scheme case 2}. Using this operator $\mathfrak{R}$ introduces another numerical error. In view of the convergence study in the previous sections that focuses on $L^1$-norm, one could aim to optimise this error measurement also in the implementation of $\mathfrak{R}$. Namely, with the notations of  Definition  \ref{de art scheme case 2} (step 2), this means that we would like to have, for each $0 \le n \le N-1$, $p \in \R^d$,
\begin{align}\label{eq error reduction}
\int |\check{v}^{N,M}_n(p,e) - \bar{v}^{N,M}_n(p,e))| \ud e
\end{align}
minimal. The solution to this problem is known, see \cite[Section 2]{jourdainreygner16}. Indeed, $\check{v}^{N,M}_n(p,\cdot)$ and  $\bar{v}^{N,M}_n(p,\cdot)$ are CDF of probability measure on $\R$, respectively denoted $\hat{\mu}^{N,M}(p)$ and $\bar{\mu}^{N,M}(p)$ for later use. It turns out that minimising \eqref{eq error reduction} corresponds exactly to minimise the $\cW_1$-distance between the two distribution.
Let us  recall the general approximation result stated in \cite[Section 2]{jourdainreygner16}. %\textcolor{red}{quote Pages?}

% Minimising \eqref{eq error reduction} amounts to find the best approximation  of $\bar{\mu}^{N,M}(p)$, which has $2dM$ atoms, by a probability distribution ($\hat{\mu}^{N,M}(p)$) which has $M$ atoms. It is shown in \cite[Section 2]{jourdainreygner16} minimising \eqref{eq error reduction} corresponds exactly to minimise the $\cW_1$-distance between the distribution.

\begin{Lemma}\label{le quantif result}
Let $F$ be the CDF associated to a probability distribution $\nu \in \cP_1(\R)$. We note $F^{-1}$ its generalized inverse, given by
\begin{align}\label{eq de gen inv}
F^{-1}(v) = \inf \set{ y \in \R| F(y) \ge v}\;,\; v \in (0,1)\,.
\end{align}
Let $\hat{\nu}_M = \frac1M\sum_{m=1}^M\delta_{e_m}$, then the minimum of
\begin{align*}
\cW_{1}(\mu,\hat{\nu}_M) = \int |H*\nu(x)-H*\hat{\nu}_M(x)|\ud x
\end{align*}
is achieved for $e_m = F^{-1}(\frac{2m-1}{2M})$, $1 \le m \le M$.
\end{Lemma}

\noindent In our framework, the above result can be specialized since the starting distribution $\bar{\mu}^{N,M}(p)$ is already atomic (with $2dM$ atoms).

%We first introduce the following definition.
%\begin{Definition}[Sorting operator] Denote $S^{M}$ the operator of sorting $M$ particles in ascending order, namely
%	\begin{align}
%		\R^{M} \ni \xi \mapsto S^{M}(\xi) \in \cD_{M}.
%	\end{align}
%\end{Definition}

\begin{Corollary}\label{eq optimising L1 norm}
For $\xi \in \cD_{2dM}$, let $\nu := \frac1{2dM}\sum_{m=1}^{2dM}\xi_m$. Then the optimal $\hat{\nu}_M$, recall Lemma \ref{le quantif result} is given by
\begin{align}
\hat{\nu}_M = \frac1M\sum_{m=1}^M \delta_{e_m} \;\text{ with }\; e_m = \xi_{(2m-1)d}\,.
\end{align}
\end{Corollary}

\proof
According to Lemma \ref{le quantif result} and the definition of $F^{-1}$ in \eqref{eq de gen inv}, one should find, for $1 \le m \le M$,
\begin{align}
i^\star(m) := \min \set{i \in \set{1,\dots,M} \, |\, F(\xi_i) \ge \frac{2m-1}{2M}}.
\end{align}
Observing that $F(\xi_i) = \frac{i}{2dM}$, we obtain that $i^\star(m) = (2m-1)d$, which concludes the proof of the corollary.
\eproof

}

\begin{Definition}[Operator $\mathfrak{R}^{M,m}$]\label{reduce particle}
	For $M\ge m \ge 1$, we first denote $S^{M}$ the operator of sorting $M$ particles in ascending order, namely
	\begin{align}
		\R^{M} \ni \xi \mapsto S^{M}(\xi) \in \cD_{M}.
	\end{align}
	At each step $n$, for a set of $2dM$ particles namely $e \in \R^{2dM}$, $S^{2dM}(e)$ could be written as follows:
	\begin{align}
		S^{2dM}(\xi) = (\widetilde S^{2d,1}(\xi),\cdots,\widetilde S^{2d,M}(\xi)),
	\end{align}
where for $1\le i\le M, \widetilde S^{2d,i}(\xi) := \Big(S^{2dM}_{j}(\xi)\Big)_{2d(i-1)< j \le 2di} \in \cD_{2d}$,  with $S^{2dM}_i(\xi)$ denoting the $i$-th coordinate of vector $S^{2dM}(\xi)$. For $\psi^\xi \in \cI^{2dM}$ associated to $\xi\in \cD_{2dM}$, we recall that 
\begin{align*}
	\psi^\xi = H\ast \Big(\frac {1}{2dM} \sum_{m=1}^{2dM}\delta_{S^{2dM}_m(\xi)}\Big),
\end{align*}
we then define the following operator
\begin{align}
	\psi^\xi \ni \cI^{2dM} \mapsto\mathfrak{R}^{2dM,M}(\psi^\xi):= H\ast \Big(\frac 1M \sum_{m=1}^{M}\delta_{\tilde e_m}\Big) \in \cI^{M},
\end{align}
where four different choices for $\tilde e_m$ are considered, leading to four different implementations. namely, for each $1\le m \le M$,
\begin{itemize}
	\item[\verb?mean?:] The mean position of $2d$ particles of $\widetilde S^{2d,m}(\xi)$, namely 
	\begin{align}
		\tilde e_m := \frac{1}{2d}\sum_{j=2d(m-1)+1}^{2dm}S^{2dM}_{j}(\xi).
	\end{align}
\item[\verb?leftmost?:]  The minimum position among the $2d$ particles of $\widetilde S^{2d,m}(\xi)$ is kept,
	\begin{align}
		\tilde e_m := \min\{S^{2dM}_{j}(\xi), 2d(m-1)+1\le j \le 2dm\} = S_{2d(m-1)+1}^{2dM}(\xi).
	\end{align}
\item[\verb?rightmost?:] The maximum position among the $2d$ particles of $\widetilde S^{2d,m}(\xi)$ is kept,
	\begin{align}
		\tilde e_m := \max\{S^{2dM}_{j}(\xi), 2d(m-1)+1\le j \le 2dm\} = S_{2dm}^{2dM}(\xi).
	\end{align}
{\color{black}
\item[\verb?optim?:] The optimal position according to Corollary \ref{eq optimising L1 norm} is kept,
	%\begin{align}
	$
		\tilde e_m := S_{d(2m-1)}^{2dM}(\xi).
	$
	%\end{align}
}	
\end{itemize}
\end{Definition}

\subsection{Numerical results}
First of all, to validate empirically our different implementations of operator $\mathfrak{R}$, we report the $L^1$-error and $L^\infty$-error against ``proxy'' solution for different numbers of particles, based on different reducing particle implementations defined in Definition \ref{reduce particle} for $\sigma=0.01, 0.3, 1.0$, see Table \ref{table case 2 sigma1}. We observe that those four implementations of $\mathfrak{R}$ including \verb?leftmost?, \verb?mean?, \verb?rightmost?, \verb?optim? of positions of particles, have practically same numerical results for different levels of volatilities but \verb?optim? is slightly better (as expected). Hence in the following, we  only consider the CASE 2 implemented with operator $\mathfrak{R}$ \verb?optim?. The clear advantage of considering the CASE 2 scheme is that the number of particles does not grow exponentially during the iterations, hence it takes less computational time compared to CASE 1, see Table \ref{ta comp time case 1 and 2}. 
% with mean position of particles.

\begin{center}
	\begin{tabular}{|c| c| c| c|c|} 
		\hline
		Sigma & Number of particles & Method & $L1$-error & $L\infty$-error \\ [0.5ex] 
		\hline\hline
		1.0&1000 &  $\verb?leftmost?$ & $5.71\mathrm{e}{-3}$ & 1.67$\mathrm{e}{-2}$\\
		\hline
		 & & $\verb?mean?$ & 4.96$\mathrm{e}{-3}$ & 1.50$\mathrm{e}{-2}$ \\
		\hline
		& & $\verb?rightmost?$ & 5.39$\mathrm{e}{-3}$&  1.80$\mathrm{e}{-2}$ \\
		\hline
		& & $\verb?optim?$ & 3.72$\mathrm{e}{-3}$&  1.02$\mathrm{e}{-2}$ \\
		\hline
		& 5000 &  $\verb?leftmost?$ & 5.05$\mathrm{e}{-3}$ & 1.45$\mathrm{e}{-2}$\\
		\hline
		& & $\verb?mean?$ & 4.92$\mathrm{e}{-3}$ & 1.46$\mathrm{e}{-2}$\\
		\hline
		& & $\verb?rightmost?$ & 4.92$\mathrm{e}{-3}$&  1.53$\mathrm{e}{-2}$\\
		\hline
		& & $\verb?optim?$ & 3.47$\mathrm{e}{-3}$&  1.01$\mathrm{e}{-2}$ \\
		\hline
		0.3&1000 &  $\verb?leftmost?$ & 1.34$\mathrm{e}{-3}$ & 0.85$\mathrm{e}{-2}$ \\
		\hline
		 & & $\verb?mean?$ & 0.88$\mathrm{e}{-3}$ & 0.61$\mathrm{e}{-2}$ \\
		\hline
		& & $\verb?rightmost?$ & 1.25$\mathrm{e}{-3}$ &  0.91$\mathrm{e}{-2}$\\
		\hline
		& & $\verb?optim?$ & 1.19$\mathrm{e}{-3}$&  0.69$\mathrm{e}{-2}$ \\
		\hline
		\hline
		& 5000 &  $\verb?leftmost?$ & 0.94$\mathrm{e}{-3}$ & 0.68$\mathrm{e}{-2}$\\
		\hline
		& & $\verb?mean?$ & 0.88$\mathrm{e}{-3}$ & 0.61$\mathrm{e}{-2}$\\
		\hline
		& & $\verb?rightmost?$ & 0.92$\mathrm{e}{-3}$ &  0.65$\mathrm{e}{-2}$\\
		\hline
		& & $\verb?optim?$ & 0.56$\mathrm{e}{-3}$&  0.46$\mathrm{e}{-2}$ \\
		\hline
		0.01 &1000 &  $\verb?leftmost?$ & 0.27$\mathrm{e}{-3}$ & 0.39$\mathrm{e}{-2}$ \\
		\hline
		 & & $\verb?mean?$ & 0.13$\mathrm{e}{-3}$ & 0.20$\mathrm{e}{-2}$ \\
		\hline
		& & $\verb?rightmost?$ & 0.48$\mathrm{e}{-3}$ &  0.43$\mathrm{e}{-2}$ \\
		\hline
		& & $\verb?optim?$ & 0.11$\mathrm{e}{-3}$&  0.17$\mathrm{e}{-2}$ \\
		\hline
		& 5000 &  $\verb?leftmost?$ & 0.09$\mathrm{e}{-3}$ & 0.13$\mathrm{e}{-2}$\\
		\hline
		& & $\verb?mean?$ & 0.08$\mathrm{e}{-3}$ & 0.17$\mathrm{e}{-2}$\\
		\hline
		& & $\verb?rightmost?$ & 0.13$\mathrm{e}{-3}$ &  0.18$\mathrm{e}{-2}$\\
		\hline
		& & $\verb?optim?$ & 0.06$\mathrm{e}{-3}$&  0.11$\mathrm{e}{-2}$ \\
		\hline
	
	\end{tabular}
	\captionof{table}{$L1$-error and $L\infty$-error for model Example \ref{ex lin model p2} with different numbers of particles with respect to \texttt{leftmost,mean,rightmost,optim}, note that the time steps $N=20$.}
\label{table case 2 sigma1}
\end{center}

%\begin{figure}
%	\begin{minipage}{.5\linewidth}
%		\centering
%		\subfloat[$\sigma=0.01$]{\label{main:a}\includegraphics[scale=.5]{img/case2_diff_sigma001.png}}
%	\end{minipage}%
%	\begin{minipage}{.5\linewidth}
%		\centering
%		\subfloat[$\sigma=0.3$]{\label{main:b}\includegraphics[scale=.5]{img/case2_diff_sigma03.png}}
%	\end{minipage}\par\medskip
%	\centering
%	\subfloat[$\sigma=1.0$]{\label{main:c}\includegraphics[scale=.5]{img/case2_diff_sigma1.png}}
%	\caption{Model of Example \ref{ex lin model p2}: Comparison of three different operators $\mathfrak{R}$ for CASE2 with $d=4$. The number of particles is $M=3500$ and the number of time steps $N= 20$.}
%	\label{fig case2 for sigma}
%\end{figure}

\begin{center}
	\begin{tabular}{|c | c| c |} 
		\hline
		Dimension & $d=1$& $d=4$\\ [0.5ex] 
		\hline\hline
		Time for CASE 1 & 46.96s  & 784s\\
	\hline
	Time for CASE 2 & 0.3s & 1.3s \\
	\hline
	\end{tabular}
	\captionof{table}{Computational cost in Example \ref{ex lin model p2} for different dimension $d$ (for the $P$-variable) for CASE 1 and CASE 2 schemes.}
\label{ta comp time case 1 and 2}
\end{center}

\begin{figure}
	\begin{minipage}{.5\linewidth}
		\centering
		\subfloat[$\sigma=0.01$]{\label{main:a}\includegraphics[scale=.5]{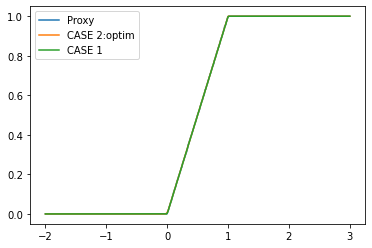}}
	\end{minipage}%
	\begin{minipage}{.5\linewidth}
		\centering
		\subfloat[$\sigma=0.3$]{\label{main:b}\includegraphics[scale=.5]{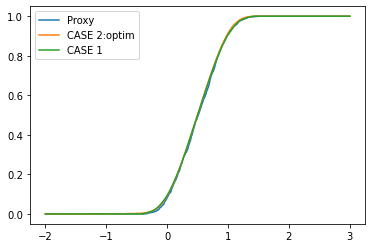}}
	\end{minipage}\par\medskip
	\centering
	\subfloat[$\sigma=1.0$]{\label{main:c}\includegraphics[scale=.5]{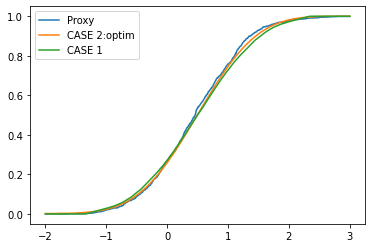}}
	\caption{Model of Example \ref{ex lin model p2}: Comparison of the two methods CASE 1 \& CASE2  with $d=4$. The Proxy solution is given by the same particle method on the one-dimensional PDE. For $BT\&SPD$: both CASE 1 and CASE 2, the number of particles is $M=3500$ and the number of time steps $N= 20$.}
	\label{fig case12 for sigma}
\end{figure}

To validate empirically our numerical schemes, we plot the value function $\cV$ of model Example \ref{ex lin model p2} using scheme CASE 1 and CASE 2 against the \emph{proxy} solution presented in \cite{bossy1996convergence}, see Figure \ref{fig case12 for sigma}. As  expected, both CASE 1 and CASE 2 scheme could reproduce correctly the entropy solution of the PDE \eqref{ass restrict regul}. We also observe that the solutions obtained by CASE 1 and CASE 2 scheme are quite close.

\begin{figure}
	\begin{minipage}{.5\linewidth}
		\centering
		\subfloat[$\sigma=0.01$]{\label{main:a}\includegraphics[scale=.5]{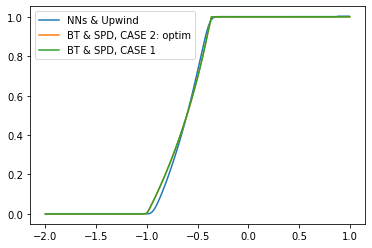}}
	\end{minipage}%
	\begin{minipage}{.5\linewidth}
		\centering
		\subfloat[$\sigma=0.3$]{\label{main:b}\includegraphics[scale=.5]{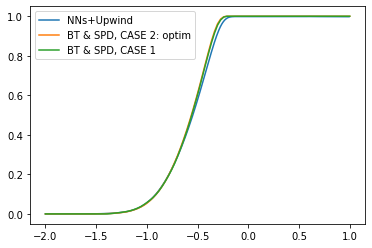}}
	\end{minipage}\par\medskip
	\centering
	\subfloat[$\sigma=1.0$]{\label{main:c}\includegraphics[scale=.5]{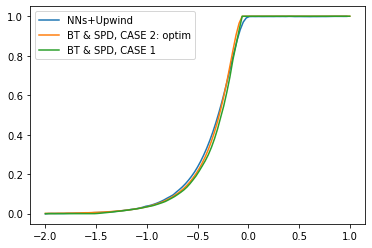}}
	\caption{Model of Example \ref{ex mulitplicative model p2}: Comparison of the two methods CASE 1 \& CASE2  with $d=4$. The Proxy solution is given by the Neural nets \& Upwind in \cite{chassagneux2022numerical}. For $BT\&SPD$: both CASE 1 and CASE 2, the number of particles is $M=3500$ and the number of time steps $N= 20$.}
	\label{fig case12 for sigma multi}
\end{figure}

Apart from the linear model, we have also tested numerically our CASE 1 and CASE 2 schemes in model Example \ref{ex mulitplicative model p2} against a NN \& Upwind method based on splitting schemes presented in \cite{chassagneux2022numerical}, see Figure \ref{fig case12 for sigma multi}. We note that the function $\mu$ in Example \ref{ex mulitplicative model p2} is always non negative, thus we use Upwind scheme which is less ``diffusive'' than Lax-Friedrisch scheme as argued in \cite{chassagneux2022numerical}. As  expected, both CASE 1 and CASE 2 schemes give satisfying numerical results for different levels of volatility.

%\textcolor{red}{TO DO: conclude this section by illustrating the rate of convergence given in the main theorem ! the empirical one should be better !!}

\begin{figure}[h]
	\centering
	\includegraphics[width=0.5\textwidth]{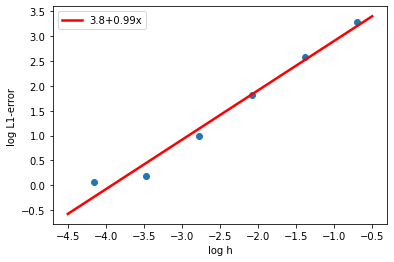}
	\caption{Convergence rate on $\mathfrak{h}$ for the model in Example \ref{ex lin model p2} with parameters $d=4,\sigma=1.0$}
	\label{theo conv btspd}
\end{figure}

At last, we want to empirically estimate the convergence rate of the error introduced by our numerical scheme. We consider the model Example \ref{ex lin model p2} where $\sigma =1.0$. We consider a set of number of time steps $N:=\set{2,4, 8, 16, 32, 64}$ and time step $\mathfrak{h}:= \frac TN$, and compute the $L1$-error for BT \& SPD method. Note that the proxy solution is always given by method in \cite{bossy1996convergence} applied to one dimensional equivalent model. The empirical convergence rate with respect to the time step is close to one see Figure \ref{theo conv btspd}, which is much better than the order $\frac 16$ obtained in Theorem \ref{th main conv result}.

\bibliographystyle{unsrt}

\bibliography{biblio}

\end{document}